\documentclass[12pt,openbib]{article}

\usepackage{amsmath,amssymb} 
\usepackage{graphicx}
\usepackage{authblk}
\newtheorem{theorem}{Theorem}
\newtheorem{lemma}{Lemma}
\newtheorem{definition}{Definition}
\newtheorem{proposition}{Proposition}[section]

\newcommand{\va}{\mathbf{a}}
\newcommand{\vb}{\mathbf{b}}
\newcommand{\vx}{\mathbf{x}}
\newcommand{\Frac}[2] {\frac{\mbox{\normalsize{$#1$}}}{\mbox{\normalsize{$#2$}}}}
\newcommand{\DP}[2]{\displaystyle \frac{\partial #1}{\partial #2}}

\newcommand{\fin}     {\hfill {\framebox[3mm]{ \rule{0cm}{1mm}}}}

\begin{document}

%
%
\title{Optimization under convexity constraints: are the finite element discretizations consistent ?}

%
%
\author[1]{Herv\'e V.J. Le Meur}
\affil[1]{\rm Laboratoire de Math\'ematiques d'Orsay, Univ. Paris-Sud, CNRS, Universit\'e Paris-Saclay, 91405 Orsay, France}



\maketitle

\begin{abstract} It is proved in \cite{PCLM} that the problem of
  minimizing a Dirichlet-like functional of the function $u_h$
  discretized with $P_1$ Finite Elements, under the constraint that
  $u_h$ be convex, cannot converge. Here, we first improve this result
  by proving that non-convergence is due to the mesh refinment lack of
  richness, remains local and is true even for any mesh. Then, we
  investigate the consistency of various natural discretizations
  ($P_1$ and $P_2$) of second order constraints (subharmonicity and
  convexity) without discussing the convergence. We also numerically
  illustrate convergence of a method proposed in the literature that
  is simpler than existing methods.
\end{abstract}

\noindent \mbox{Keywords:convexity; finite elements; interpolation; conformal approximation; minimization}


\noindent \mbox{subjclass:{26B25, 52A41, 65K10, 65N30, 90C25, 91B24}}
\section{INTRODUCTION}

This paper is devoted to the numerical discretization of optimization
with constraints on the second order derivative, namely problems of
the form

\begin{equation}
\label{eqm4}
\left\{
\begin{array}{l}
\inf J(u),\\
\mathbf{D}^2 u\in K,
\end{array}
\right.
\end{equation}

\noindent where $J$ is a functional and $K$ is a subset of the set of
symmetric matrices. Such problems appear in various contexts, in
particular in physics and economics.

\bigskip

Since Newton, the shape of minimal resistance has been a topic of
interest. It is called the Newton's problem and is of the type of
(\ref{eqm4}). With some additional assumptions, it amounts to looking
for a concave function that minimizes a nonlinear functional. See for
instance the original book \cite{Newton}, the historical survey in
\cite{Goldstine}, more recent theoretical results in
\cite{buttazo,LR_Pelletier_01,PL2} or numerical results
\cite{LR_Oudet_05,Waschumth_14}.

The associated problem of discretizing convex functions or bodies
is indeed wider than expected. For instance, the Alexandrov's
problem (see \cite{Carlier_04} and \cite{LR_Oudet_05}), the Cheeger's
constant \cite{LR_Oudet_05,Carlier_Comte_Peyre_09} and the Newton's problem
\cite{LR_Oudet_05,Waschumth_14} can be numerically studied.\\

In economics, it suffices to remember that utility functions are
concave to see the very wide applicability of these problems. For
instance, in \cite{RC}, the authors are interested in finding the
minimum of a convex and quadratic functional $J$ over the set of
convex functions
\begin{equation}
\label{eqm3}
\min_{u \in K} J(u) \mbox{ for } \; J(u)=\int_{\Omega}\left( \Frac{1}{2}\mathbf{\nabla} u^T \, \mathbf{C} \, \mathbf{\nabla} u -\mathbf{x} \cdot \mathbf{\nabla}
u+(1-\alpha)u\right){\rm d} \mathbf{x},
\end{equation}
where $0\leq\alpha\leq 1$ and $ \mathbf{C}$ is a (2,2) symmetric and
positive-definite matrix and $K$ is given by
\begin{equation}
\label{eqm2}
K=\{u\in H^1(\Omega), u\geq 0,u_x\geq 0,u_y\geq 0, u \mbox{ convex }\}.
\end{equation}
So $J$ is strictly convex and the set $K$ is convex. It is then easy to
check that there exists a unique minimizer of $J$ over $K$. The
regularity of the solutions is studied in \cite{Carlier_Lachand-Robert_01}.
Note that, when $\alpha=1$ and in a domain $\Omega$ with nonnegative
coordinates, the problem degenerates: an exact convex solution exists,
up to an additive constant
\begin{equation}
\label{eqm1}
\mathbf{\nabla} u(\mathbf{x})=\mathbf{C}^{-1}\mathbf{x} \Rightarrow u(\mathbf{x})= \mathbf{x}' \mathbf{C}^{-1} \mathbf{x}/2 -\mbox{Cst},
\end{equation}
that can be fixed if we force $u$ to vanish at a point (if
$\alpha=1$).

In \cite{PCLM}, Chon\'e and Le Meur exhibit an obstruction to the
convergence of the mere {\em discretization} of this problem through
conformal $P_1$ Finite Element (FE) method on regular meshes. They
also extensively use the explicit solution (\ref{eqm1}) to illustrate
the lack of convergence when a predicted condition is no more
satisfied. Here ``conformal'' means that the {\em discretized}
function is supposed to satisfy exactly the continuous constraint. As
a consequence of this result, no optimization process that would use
such conformal $P_1$ discretization of both the functional and the
constraint may converge for any solution on these meshes. Yet
approximation theory easily proves convergence of every separate
discretization of these terms.\\

A possible approach to circumvent the mesh problem is, first, to test
whether a sample of values on {\em given points} (and not {\em mesh})
may be associated to a convex function or body. Then one must
construct the associated mesh (so, function dependent !) so as to
interpolate. This was done in \cite{Leung_Renka_99} for $C^1$ FE by
Leung and Renka. But such a regularity is too restrictive for us. This
article \cite{Leung_Renka_99} reviews other papers, some of them being
commented as false. It proves that the issue is not so simple.

The $H_0^1$ projection of a given function is addressed in
\cite{Carlier_LR_Maury_99} through the saddle-point method. In this
article, the authors even weight the convexity constraint to tune the
convergence to the solution. The computations appear to be more robust
than expected by theory. An attempt of explanation is given in
\cite{Maury_03}.

In \cite{Carlier_LR_Maury_01}, the authors describe and implement an
algorithm that optimizes not in the set of discretized {\em and}
convex (i.e. conformal) functions but in the set of the convex
functions after discretization (so that they may be non-convex once discretized). More
precisely they characterize, for a specific structured set of
cartesian points, the image through the $P_1$ discretization of a
(continuous) convex function. This yields such a huge number of
constraints on the function that it may not be recommended.

In \cite{LR_Oudet_05}, Lachand-Robert and Oudet address the problem of
discretizing a convex body. They use the parametric generalization in order to
discretize convex functions' graphs. In the functional, they isolate
the dependence on the point $x$, the unit normal $\nu$ at $x$ and the
signed distance $\phi= \nu . x$. Although they notice that the three
variables are ``somehow redundant'', they implement a gradient like
method, based on the variation of only $\phi= x. \nu$. They address
both the Alexandrov's problem, the Cheeger's sets and the Newton's
problem. This was revisited since in \cite{Carlier_Comte_Peyre_09} and solved in a very elegant way in \cite{Waschumth_14}.

Separately, Aguilera and Morin \cite{Aguilera_Morin_08} prove
convergence of a Finite Difference (FD) ``approximation using positive
semidefinite programs and discrete Hessians''. In
\cite{Aguilera_Morin_09}, the same authors prove even convergence of a
Finite Element (FE) discretization of the weak Hessian. Since FD are
included in FE for convenient meshes, these two articles seem to
contradict both \cite{PCLM} and the present article. Indeed, they do
not, but we discuss them below (section \ref{sec.5}).

In \cite{Ekeland_Moreno_10}, Ekeland and Moreno propose a non-local
discretization that relies on the representation of any convex
function as the supremum of affine functions (its minorants). So their
``discrete'' representation of continuous functions is conformal but
non-local. The major drawback is that the complexity increases
drastically with dimension, due to the nonlocality, but it works.

More recently, two very different means of resolution were
proposed. Both are non-conformal in the sense that the discretized
solution is not convex. In the first one, M\'erigot and Oudet
\cite{Merigot_Oudet_14} propose to discretize the convexity constraint
by under sampling it. More precisely, they prove convergence of their
algorithm if the constraint is forced only on a subset of the sampling
points. In the second one, Mirebeau \cite{Mirebeau_14} proves that for
a given grid, he may provide a sequence of sets of (less than four)
points on which it suffices to force convexity to ensure convexity at
convergence on this given grid. This makes only $O(N\log^2 N)$ linear
constraints.

There seems to be a convergent opinion through \cite{Mirebeau_14},
\cite{Merigot_Oudet_14}, \cite{Aguilera_Morin_09} that indeed the
number of constraints must be not too large. A deeper discussion of
this crucial point is postponed to a forthcoming article.

In \cite{Waschmuth_17}, the author goes further than
\cite{Aguilera_Morin_09}. He uses approximation theory, where there is
no reason why a function (should it be convex or not) should not be
approximable. So as to remain inside this theory, he uses two
auxilliary functions that make his solution strictly convex, so that
the constraint of convexity may not be saturated. Nevertheless, we are
surprised that these functions are not needed in his numerical
experiment. It could be due to the fact that the hidden constraints
(that will appear below) are filled by the exact solution and so do
not trigger any discrepancy. This explains why our results are not
contradictory with Oberman's or Wachsmuth's.

A totally different tool has emerged since. It consists in using an
evolution PDE with a $u_0$ as an initial condition. In
\cite{Oberman_08} (and other articles), the author computes the convex
envelope of a graph using such a nonlinear PDE. Later
\cite{Carlier_Galichon_11} exhibits an evolution PDE, close to the
previous, which, starting from a given $u_0$, is proved to converge to
the convexified of $u_0$. They use stochastic control representation
to prove an exponential convergence. This set of evolution PDE is a
new idea and very important in optimal transport theory.

The present article deals with two main types of non-convergence
results. On the one hand, the non-convergence in $L^2$ for conformal
$P_1$ Finite Element (FE) is revisited after \cite{PCLM} from a
theoretical point of view. On the other hand, the consistency of
discretized linear optimization over second order constraints is
investigated for various discretizations ($P_1$ and $P_2$) and various
constraints: linear (subharmonicity) or nonlinear (convexity). The
motivation for convexity has been stressed. Subharmonicity is only
studied as an intermediate step toward convexity because it is
linear.\\

In Section \ref{sec.2}, first we restate some already known results,
then we prove that non-convergence is purely local and so is true for
any mesh. In Section \ref{sec.3}, we investigate the consistency of
various $P_1$ FE discretizations of our model problem. Section
\ref{sec.4} is devoted to the consistency of $P_2$ FE discretizations
(strong or weak convexity). We discuss the articles
\cite{Aguilera_Morin_08,Aguilera_Morin_09} that could seem to
contradict \cite{PCLM} in Section \ref{sec.5} and we conclude in
Section \ref{sec.6}.

\section{THEORETICAL RESULTS ON THE $P_1$ FEM}
\label{sec.2}

In the present section, we first recall some already known results and
make them more explicit.

\subsection{Already known results}


In \cite{PCLM}, various results are proved. Since the goal in the
present subsection is to extend this study, we need to remind the
reader of these results. We consider an open convex and bounded domain
$\Omega$ of $\mathbb{R}^2$, but we use only $\Omega = (a,b)^2$ (and
$a$ and $b$ positive) in applications. The Lemma 3 of \cite{PCLM}
states:
\begin{lemma}
\label{lem1}
A function $u_h$, $P_1$ in every triangle of $\Omega$'s mesh, is
convex if and only if, for any pair of adjacent triangles
\begin{equation}
\label{eq1}
(\mathbf{q}_2 -\mathbf{q}_1) \cdot \mathbf{n} \geq 0,
\end{equation}
where $\mathbf{q}_1$ (resp. $\mathbf{q}_2$) is the (constant) gradient
of $u_h$ inside triangle 1 (resp. 2) and $\mathbf{n}$ is the unit
normal pointing from triangle 1 to 2.
\end{lemma}

As reminded in \cite{PCLM}, convexity of $u$ can be defined dually as
$-\langle \mathbf{\nabla}u \otimes \mathbf{\nabla}v\rangle \succeq 0,$
where $\mathbf{A} \succeq 0$ means the matrix $\mathbf{A}$ is positive
semi-definite, $v$ is a nonnegative test function in
$\mathcal{C}_0^{\infty}$ and
\begin{equation}
\label{def_otimes}
-\langle \mathbf{\nabla}u \otimes \mathbf{\nabla}v\rangle  =-\left( \begin{array}{cc} \langle u_{x} ,v_{x} \rangle & \langle u_{x} , v_{y} \rangle\\ \langle u_{y} , v_{x} \rangle &  \langle u_{y} , v_{y} \rangle \end{array} \right) .
\end{equation}
This definition is weak (in the sense of distributions) but it is also
equivalent to a strong one when $u$ is sufficiently regular. We will
use a weaker definition where the trial and test functions will be in
a (finite dimensionnal) subspace of $\mathcal{C}_0^{\infty}$ : the set
of some FE functions. Last we will report on \cite{Aguilera_Morin_09}
which uses an even weaker definition where the test functions are
in a much smaller subspace than $u_h$'s. One can easily check that the
weak definition of convexity (with ${\mathcal C}_0^{\infty}$ test
functions) implies the weaker definition of convexity (with only the
basis functions in $P_1$ associated to interior points as test
functions) but they are not equivalent.

The proof of Lemma \ref{lem1}  is based on the following formula written
for the $P_1$ function $u_h$:
\begin{equation}
\label{eq2}
\langle \frac{\partial^2 u_h}{\partial \va \partial \vb},\varphi\rangle =\sum_e (\mathbf{q}_2-\mathbf{q}_1) \cdot \mathbf{n} \, (\mathbf{n} \cdot \va)(\mathbf{n} \cdot \vb)\int_e\varphi(s){\rm d} s,
\end{equation}
where the summation is taken over all interior edges $e$ of the mesh,
$\va$ and $\vb$ are two unit vectors, $\varphi$ is a ${C}^{\infty}$
function with compact support in $\Omega$, $\mathbf{q}_1$ and
$\mathbf{q}_2$ designate the two (constant) gradients of $u_h$ in the two
triangles that share the edge $e$ and $\mathbf{n}$ is the unit normal
from triangle 1 to triangle 2. This formula is right thanks to the
property that the gradient of a continuous function across the edge
has its tangential derivative along the edge continuous. Namely
$(\mathbf{q}_2-\mathbf{q}_1) \cdot \mathbf{t}=0$ where $\mathbf{t}$ is the
unit tangent vector. The proof also relies on the fact that for all unit vector
$\mathbf{n}$, the bilinear form $(\va, \vb) \mapsto
(\mathbf{n} \cdot \va)(\mathbf{n}\cdot \vb)$ is positive semi-definite, whatever
$\mathbf{n}$.\\

Through a proof very similar to the one of Theorem 4 in \cite{PCLM},
one may prove a wider Theorem that even applies locally:
\begin{theorem}
\label{th2}
Let $\Omega$ be an open and convex subset of $\mathbb{R}^2$ and $\mathcal{T}_h$ a
triangulation of $\Omega$. If there exists an open and convex $\Omega'
\subset \Omega$ such that the following property is satisfied in
$\Omega'$ and for the mesh $\mathcal{T}_h$:
\begin{equation}
\label{def_PM}
(PM)
\left\{
\begin{array}{l}
\exists (\va ,\vb ) \mbox{ two independent unit vectors such that } (\mathbf{n} \cdot  \va )(\mathbf{n} \cdot  \vb )\geq 0 \\
\mbox{ for all } \mathbf{n} \mbox{ unit normal to an edge of } \mathcal{T}_h \bigcap \Omega',
\end{array}
\right.
\end{equation}
then any function
$u_h$ convex and $P_1$ will satisfy the following
equation in the sense of distributions on $\Omega'$:
\begin{equation}
\label{eq4}
\frac{\partial^2 u_h }{\partial \va \partial \vb } \geq 0.
\end{equation}
\end{theorem}

The proof relies on the fact that if a $P_1$ function $u_h$ is convex, then by
Lemma \ref{lem1}, the gradient's jumps are non-negative. As a
consequence, in Formula (\ref{eq2}), if $u_h$ is convex, the scalar
coefficients $(\mathbf{q}_2-\mathbf{q}_1)\cdot \mathbf{n}$ are all
non-negative. But separately, if condition $(PM)$ is satisfied,
$(\mathbf{n}\cdot \va)(\mathbf{n}\cdot \vb)\geq 0$ and so even the non-diagonal
terms of the hessian ($\partial^2 u_h / (\partial \va \partial \vb)$)
are sign-constrained !

Notice first that property $(PM)$ only depends on the geometry of the
mesh and not on any function. Notice then that the property $(PM)$ is
at given $h$ but, if it remains for a sequence of $h \rightarrow 0$ in
$\Omega' \subset \Omega$, then the associated functions $u_h$ given by
Theorem \ref{th2} satisfy (\ref{eq4}) even at the limit in $\Omega'
\subset \Omega$. Yet such a property (\ref{eq4}) is (sometimes)
contradictory with the property of convexity. For instance $(x,y)
\mapsto (x^2-2\rho x y+\mu y^2)/(\mu -\rho^2)$, where $\mu-\rho^2>0$
and $1+\mu>0$, is strictly convex. Yet $\partial^2 u/(\partial
x \partial y)=-\rho/(\mu-\rho^2)$ has no sign prescribed if $u$ is
only strictly convex. This is the key of the obstruction to the
convergence stated in \cite{PCLM}. The next subsection is devoted to
discussing the generality of $(PM)$.

\subsection{Is property $(PM)$ frequent ?}

\begin{figure}[hbtp]
    \centering
    \includegraphics[height=4cm]{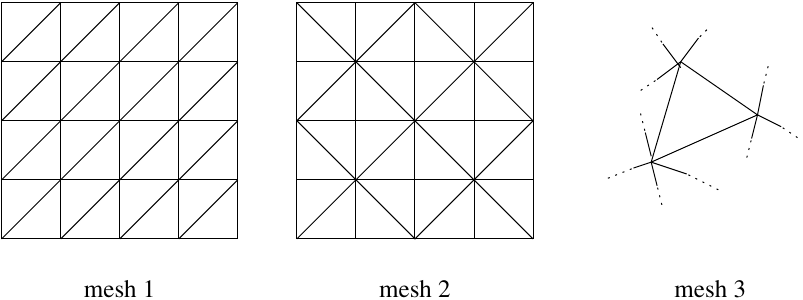}
    \caption{Three meshes.}
  \label{fig1}
\end{figure}

In this subsection, we start from some elementary observations on
particular meshes, then we argue on whether such meshes are frequent
and state a Theorem. Three types of meshes may be considered that are
depicted in Figure \ref{fig1}.


Concerning mesh 1, three different unit normals (up to a
multiplicative factor $-1$) exist in {\em all} the domain: $\{(0,1);
(1,0) ; (-1/\sqrt{2}, 1/\sqrt{2}) \}$. If we choose $\va = (-1,0)$ and
$\vb = (0,1)$, the values taken by $(\mathbf{n} \cdot \va )(\mathbf{n}
\cdot \vb )$ (with the possible $\mathbf{n}$) are $\{0 ; 0 ;
(1/\sqrt{2})(1/\sqrt{2}) \} $. They are all nonnegative.

Concerning mesh 2, four types of unit normals can be found in {\em
  all} the domain: $\{(0,1);$ $ (1,0) ;$ $ (-1/\sqrt{2}, 1/\sqrt{2})
;$ $ (1/\sqrt{2}, 1/\sqrt{2})\}$. If we choose $\va = (1,0)$ and $\vb
=(1/\sqrt{2}, -1/\sqrt{2})$, the values taken by $(\mathbf{n} \cdot
\va )(\mathbf{n} \cdot \vb )$ (with the possible $\mathbf{n}$) are
$\{1/\sqrt{2} ; 0 ; (-1/\sqrt{2})(-1) ; 0 \} $. They are all
nonnegative.

Any structured mesh like mesh 1 or 2, if refined while
keeping the same structure, will satisfy $(PM)$ with the very same
$\va$ and $\vb$ for any $h$. But what can be stated about a more
general mesh like mesh 3 ?

We are going to prove the following Theorem.

\begin{theorem}
\label{th3}
Let $\Omega$ be an open and convex set of $\mathbb{R}^2$, $\Omega' \subset
\Omega$ an open and convex subset and $\mathcal{T}_h$ a triangulation of
$\Omega$.\\
For any given $h$ and $\Omega'$, there exists $(\va ,\vb )$ such that
$(PM)$ is true in $\Omega'$ and for $\mathcal{T}_h$. Moreover, if the
refinement process does not enrich the edges' directions of
$\mathcal{T}_{h_n}$ in $\Omega'$, then the property $(PM)$ is
true for any $\mathcal{T}_{h_n}$ in $\Omega'$.
\end{theorem}

\noindent{\em Proof of Theorem \ref{th3}}

\noindent First, one must notice that the condition $(\mathbf{n} \cdot \va
)(\mathbf{n} \cdot \vb ) \geq 0$ is invariant through changing
$\mathbf{n}$ into $-\mathbf{n}$. So, it does not depend on the choice of
the hyperplane's unit normal.

\begin{figure}[hbtp]
    \centering
    \includegraphics[height=6cm]{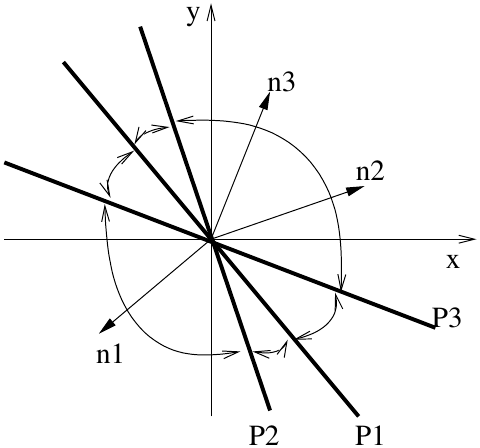}
    \caption{Unit normals and the associated hyperplanes in case of a triangle.}
  \label{fig2}
\end{figure}

Let us take a very general mesh like mesh 3 for which we may isolate a
single triangle. Then, for each of the three unit normals (up to the
multiplicative coefficient -1), having $(\mathbf{n} \cdot \va )(\mathbf{n}
\cdot \vb ) \geq 0$ is equivalent to having $\va$ and $\vb$ in the {\em
  same} closed half plane whose normal is $\mathbf{n}$. Equivalently,
we can take either $\va$ and $\vb$ on the same side as $\mathbf{n}$ or
on the opposite side. Gathering the conditions associated to each of
the three edges, we are led to choosing both $\va$ and $\vb$ in one of
the six half cone that intersects none of the three hyperplanes. Such a
configuration is depicted on Figure (\ref{fig2}). The three unit
normals are drawn and called n$1$, n$2$, n$3$ and their associated
hyperplanes are called P$1$, P$2$, P$3$. So, one may easily
find such a couple $(\va, \vb)$ for any given triangle.

Even better, whatever might be the finite number of edges for a larger
$\Omega'$ and a given $\mathcal{T}_h$, one sees obviously that such a
choice of $\va$ and $\vb$ is easy since there always exists a finite
number of hyperplanes. As a consequence, $\va$ and $\vb$ can be
choosen to be {\em both} in any of the cones which partition the whole
unit circle. So the property $(PM)$ is true even on any {\em given}
and potentially unstructured mesh of a subdomain of $\Omega$.
%

Could a refinement process preserve the property $(PM)$ for a given couple
($\va, \vb$), a given $\Omega'$ and any $h \rightarrow 0$ ?

In case of any general given mesh $\mathcal{T}_h$ we have just proved
that one may find $\va$ and $\vb$, should they suit the property only
in $\Omega'$. Then, for instance, if the refinement process is such
that every triangle in $\Omega'$ is refined into four homothetic
subtriangles, it is obvious that no more edges' direction will be
provided. As a consequence, the very same $\va$ and $\vb$ will then
suit property $(PM)$ in $\Omega'$ for {\em any} refined mesh and {\em
  any} $h$. The proof is complete.

\fin


In a sense to be defined, the set of refinement processes that enable
$(PM)$ is open and non-empty. The main question is ``Does the
refinement process enrich in normals $\mathbf{n}$ ?''.


%

\subsection{More information on the lack of convergence}

\label{section2.1.3}

In \cite{PCLM}, the authors state that the ''conformal method may not
converge for some limit function'' because the second derivative of
the limit is forced to satisfy an unnatural condition. In this
subsection, we give a more precise result in $L^2$ and still use our
Lemma \ref{lem1} that enables to identify the convexity of $u_h$ and
the non-negativity of its gradients' jumps (for $P_1$ FE). This
enables us to state our main Theorem:

\begin{theorem}
\label{th6}
Let $\Omega$ be an open and convex domain in $\mathbb{R}^2$, $\Omega'$ an open and convex
subdomain of $\Omega$ and a family of meshes $(\mathcal{T}_h)_{h
 \rightarrow 0^+}$. If the property $(PM)$ is satisfied in $\Omega'$ and for the meshes $(\mathcal{T}_h)_{h \rightarrow 0}$, then there exists $\varepsilon
>0$ and a $\mathcal{C}^{\infty}$ convex function $u_{exact}$ such that
\[ 
\min_{u_h  \in P_1 \texttt{ and } u_h \texttt{ convex}} \mid u_h-u_{exact} \mid_{L^2} \geq \varepsilon.
\]
\end{theorem}

This Theorem uses that if Property $(PM)$ is satisfied in a domain
$\Omega$ for a family of meshes $(\mathcal{T}_h)_{h \rightarrow 0^+}$,
should it be only locally in $\Omega'$, there will be two unit vectors
$\va, \vb$ such that (\ref{eq4}) holds. Given those $\va, \vb$, there
exists a function $u_{exact}$ that may not be the limit in $L^2$ of
{\em any} sequence of functions both $P_1$ and convex (in a strong
definition) on {\em this} family of meshes. The improvement of this
Theorem with respect to \cite{PCLM} is that the obstruction is
local. As a consequence we may apply our obstruction on
{\em any} mesh and not only regular ones.

We need some more preliminary results and definitions before proving
this Theorem \ref{th6} concerning $P_1$ FE.

\begin{definition}
\label{def1}
Let $\Omega$ be an open subset of $\mathbb{R}^2$, and $u \in
\mathcal{C}^0(\Omega)$. For any $\vx = (x,y) \in \Omega$ and $\{\va ,
\vb \}$ two independent unit vectors, there exist two positive numbers
$\{\alpha_0,\beta_0 \}$ such that
\[
\begin{array}{c}
\forall \alpha , \;  0 \leq \alpha \leq \alpha_0, \forall \beta , \; 0\leq \beta \leq \beta_0, \\
\vx+\alpha \va+\beta \vb = (x+\alpha a_1+\beta b_1,y+\alpha a_2+\beta b_2) \in \Omega.
\end{array}
\]
For such $(\alpha_0, \beta_0, \va, \vb )$, one defines:
\[
\phi_{(\alpha_0, \beta_0, \va, \vb )}(u)=\left(u(\vx+\alpha_0 \va +\beta_0 \vb)-u(\vx+\alpha_0 \va )-u(\vx +\beta_0 \vb)+u(\vx)\right)/(\alpha_0 \beta_0).
\]
\end{definition}

If $\alpha_0 = \beta_0 \rightarrow 0^+$, then
$\phi_{(\alpha_0, \beta_0, \va, \vb )}(u) \rightarrow \partial^2 u
/(\partial \va \partial \vb)$. We are going to prove that such a quantity
$\phi_{(\alpha_0, \beta_0, \va, \vb )}(u)$ will be overconstrained by
the mere discretization and the limit $u$ will not satisfy the right
sign of the second derivative.\\

We want now to define an explicit solution depending on some
parameters. A well-chosen combination of these parameters will provide
a non-approximable function. For instance, let $\Omega=(1,2)^2$ and the problem
(\ref{eqm3}, \ref{eqm2}) for a positive definite real matrix $\mathbf{C}$
be such that:

\begin{equation}
\label{eq5}
\mathbf{C} = \left( \begin{array}{cc} \mu_2 & \rho \\ \rho & \mu_1 \end{array}\right).
\end{equation}
Then there exists an exact solution thanks to (\ref{eqm1}):
\begin{equation}
\label{eq5.5}
\begin{array}{rcl}
u_{exact} & = & \Frac{1}{\mu_1\mu_2 -\rho^2}\left(\mu_1 (x^2-1)/2+\mu_2(y^2-1)/2-\rho (x y -1)\right) \mbox{ in } \Omega.
\end{array}
\end{equation}
The function $u_{exact}$ is such that it is zero at the corner $(1,1)$
of the domain $\Omega$ chosen but it can easily be generalized to
other domains. Simple computations prove the following formula:
\begin{equation}
\label{eq6}
\phi_{(\alpha_0, \beta_0, \va, \vb )}(u_{exact})=\Frac{1}{\mu_1 \mu_2 -\rho^2}\left(\mu_1 a_1 b_1 +\mu_2 a_2 b_2 -\rho (a_1 b_2 +a_2 b_1)\right) = \va' \mathbf{C}^{-1} \vb,
\end{equation}
where $a_1, a_2, b_1, b_2$ are the components of $\va, \vb$. It is
then useful to state the following Lemma.

\begin{lemma}
\label{lem2}
Let $\va, \vb$ be two given independent unit vectors in $\mathbb{R}^2$. Let
$\Omega \subset \mathbb{R}^2$ and $\eta>0$ given. Then, there exists $(\mu_1,\mu_2,\rho)$ such that
$\mu_1\geq 0, \mu_2\geq 0$, $\mu_1 \mu_2 -\rho^2\geq 0$ and
\begin{equation}
\label{eq7}
\phi_{(\alpha_0, \beta_0, \va, \vb )}(u_{exact}) = \va^T \mathbf{C}^{-1} \vb \leq -\eta \mbox{ in }\Omega,
\end{equation}
provided $(\alpha_0,\beta_0) \in \mathbb{R}^{+*2}$ are such that there exists $\phi_{(\alpha_0, \beta_0,   \va, \vb )}(u_{exact})$.
\end{lemma}

Roughly speaking, Lemma \ref{lem2} states that for any independent
$\va, \vb$, one may find a positive definite matrix $\mathbf{C}^{-1}$
such that the associated $u_{exact}$ satisfies the following property
on its associated bilinear form:
\[
\phi_{(\alpha_0, \beta_0, \va, \vb)}(u_{exact})=\va^T \mathbf{C}^{-1} \vb <0.
\]

\noindent {\em Proof of Lemma \ref{lem2}}

\noindent From the formula (\ref{eq6}), one sees that it is sufficient
to find a positive definite symmetric bilinear form
($\mathbf{C}^{-1}$) such that, for independent $\va, \vb$ given, $\va^T
\mathbf{C}^{-1} \vb <0$. We will build up $\mathbf{C}^{-1}$ from its
eigenvalues and eigenvectors.

Let $\mathbf{e}_1$ be a unit vector between $\va$ and $\vb$
(normalized mean for instance). Then let $\mathbf{e}_2$ be a unit
vector normal to $\mathbf{e}_1$. With the same notations for $\va$ and
$\vb$ in this new basis, $a_1 b_1 > 0$ and $a_2 b_2 <0$.

Let now $\mathbf{C}^{-1}$ be a matrix in the canonical basis with the
eigenvectors $\mathbf{e}_1, \mathbf{e}_2$ and the associated positive
eigenvalues $\lambda_1, \lambda_2$. Obviously, $\mathbf{C}$ is
positive semi-definite. Then $\va^T \mathbf{C}^{-1} \vb =\lambda_1 a_1
b_1 + \lambda_2 a_2 b_2$ which may be less than $-\eta$ given, for an
appropriate choice of $\lambda_1, \lambda_2$. Then the $(\mu_1, \mu_2,
\rho)$ are the coefficients of the matrix $\mathbf{C}$ in the
canonical basis as written in (\ref{eq5}).

\fin

We have now obtained all what is needed to start the following proof.\\

\noindent {\em Proof of Theorem \ref{th6}}.

\noindent Since we assume $(PM)$ is satisfied in $\Omega'$ for
$(\mathcal{T}_h)_{h \rightarrow 0}$, this provides a couple of vectors
$(\va ,\vb )$. Theorem \ref{th2} enables us to claim that any $P_1$
function $u_h$ satisfies (\ref{eq4}) in the sense of distributions in
$\Omega'$:
\begin{equation}
\label{eq8}
0 \; \leq \; \langle \Frac{\partial^2 u_h}{\partial \va \partial \vb}, \Psi\rangle _{\Omega}=\langle u_h,\Frac{\partial^2 \Psi}{\partial \va \partial \vb}\rangle _{\Omega},
\end{equation}
for any nonnegative $\Psi \in \mathcal{C}^{\infty}_0(\Omega')$. This
property obviously remains for open and convex subdomains of $\Omega'$.

Should it be needed, one could decrease $\alpha_0, \beta_0 >0$ and
find a subdomain $\Omega'' \subset \Omega'$ such that for any
nonnegative $\varphi \in \mathcal{C}^{\infty}_0(\Omega'')$, the
function
\begin{equation}
\label{eq9}
\Psi : \; \vx \mapsto \int_0^1 \int_0^1 \varphi(\vx -t \alpha_0 \va -t' \beta_0 \vb) \, {\rm d}t \, {\rm d}t',
\end{equation}
is nonnegative and in $\mathcal{C}^{\infty}_0(\Omega')$. So the
function $\Psi$ is eligible for equation (\ref{eq8}) whatever $\phi$. The second
derivative of $\Psi$ may be computed:
\[
\begin{array}{rl}
\Frac{\partial^2 \Psi}{\partial \va \partial \vb} & = \left( \varphi(\vx -\alpha_0 \va -\beta_0 \vb)-\varphi(\vx -\beta_0 \vb)-\varphi(\vx -\alpha_0 \va)+\varphi(\vx)\right)/(\alpha_0 \beta_0)\\
 & =\phi_{(\alpha_0, \beta_0, -\va, -\vb )}(\varphi).
\end{array}
\]
So we have, for any nonnegative $\varphi \in
\mathcal{C}^{\infty}_0(\Omega'')$:
\begin{equation}
\label{eq10}
\begin{array}{rl}
0 \leq \langle \Frac{\partial^2 u_h}{\partial \va \partial \vb}, \Psi\rangle _{\Omega}& =\langle u_h,\Frac{\partial^2 \Psi}{\partial \va \partial \vb}\rangle _{\Omega}\\
 & =\langle u_h, \phi_{(\alpha_0, \beta_0, -\va, -\vb )}(\varphi)\rangle _{\Omega'}\\
 & = \langle \phi_{(\alpha_0, \beta_0, \va, \vb )}(u_h),\varphi\rangle _{\Omega'}.
\end{array}
\end{equation}

Separately, given $\va, \vb, \eta > 0$, Lemma \ref{lem2} enables us to
claim there exists a convex quadratic function $u_{exact}$ such that
(\ref{eq7}) is true for $\eta >0$ given. Moreover, because of
(\ref{eq10}), since the open set $\Omega'' \subset \Omega'$ is of
non-zero measure, and $\phi_{(\alpha_0, \beta_0, \va, \vb
  )}(u_{exact})$ is a constant known by (\ref{eq6}), we have
\[
\langle \phi_{(\alpha_0, \beta_0, \va, \vb
)}(u_h-u_{exact}),\varphi\rangle _{\Omega''} \; \geq -\phi_{(\alpha_0, \beta_0,
\va, \vb )}(u_{exact})\int_{\Omega} \varphi.
\]
As a consequence, for convenient $\Omega''$, we are led to
\begin{equation}
\label{eq12}
\langle \phi_{(\alpha_0, \beta_0,
\va, \vb )}(u_h-u_{exact}),\varphi\rangle _{\Omega''} \geq \eta \int_{\Omega} \varphi,
\end{equation}
for any {\em nonnegative} $\varphi \in
\mathcal{C}_0^{\infty}(\Omega'')$ and $u_h$ convex and $P_1$. This
last inequality contradicts any possible convergence in
$L^2(\Omega'')$ of a sequence $u_h$ of convex functions $P_1$ in
$(\mathcal{T}_h)_h$ to $u_{exact}$ given by Lemma \ref{lem2}.

\fin

The previous Theorem adds one more argument to the need of convenient
numerical methods to discretize the constraint of convexity.

\section{Consistency of the $P_1$ FEM}
\label{sec.3}

The non-convergence of the convexity problem (\ref{eqm3}, \ref{eqm2})
discretized by conformal $P_1$ Finite Elements is proved and
numerically illustrated in \cite{PCLM}. It is even proved to apply to
non-regular meshes in the subsection \ref{section2.1.3}. We give here
a very different argument by studying the consistency of various
discretizations of convexity.

One could wonder why study consistency. The reason is that when one
discretizes convexity (or even its linear form of subharmonicity),
there exist meshes and FE on which the discretization is not
consistent. As a consequence, any proof of convergence must exclude at
least these cases. The fact that the FE-discretized Dirichlet
functional and the FE-discretized Hessian matrix converge separately is not
sufficient. Already \cite{KS} claimed convergence, but was later
proved to be wrong in \cite{PCLM}.

\subsection{Use of the gradients' jumps for convexity}

\label{sec3.1}

\begin{figure}[hbtp]
    \centering
    \includegraphics[height=6cm]{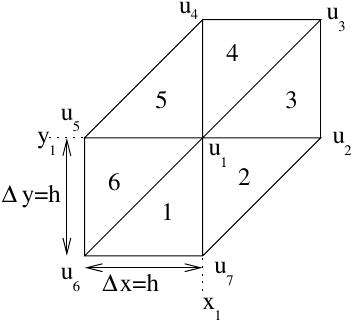}
    \caption{Local shape of the mesh close to $(x_1,y_1)$.}
  \label{fig3}
\end{figure}

Given a sufficiently smooth function $u$, one may interpolate it in
$P_1$ FE as $u_h= \sum_{i=1}^N u_{hi} \phi_i(x)$. If the mesh is
structured as in Figure \ref{fig3}, should it be local, one may
compute the gradients' jumps accross the edges as functions of the
values at the various involved nodes. Then one may compute the series
expansion of the sampled values $u_i$ from the exact initial function
$u$. The jumps between triangles 1 and 2, 2 and 3 and last 3 and 4 are
respectively:
\begin{equation}
\label{eq13}
\begin{array}{rrl}
\mbox{Jump(1/2)} =&\frac{+u_6-u_7-u_1+u_2}{h} & = (u_{xx}+u_{xy})_{x=x_1,y=y_1}h +O(h^2),\\
\mbox{Jump(2/3)} =&\frac{-u_1+u_7-u_2+u_3}{h} & = (u_{xy}+u_{yy})_{x=x_1,y=y_1}h +O(h^2), \\
\mbox{Jump(3/4)} =&\sqrt{2}\frac{-u_1+u_2-u_3+u_4}{h} & = -\sqrt{2}(u_{xy})_{x=x_1,y=y_1}h +O(h^2),
\end{array}
\end{equation}
where we denote $u_x$ the partial derivative with respect to $x$ and
$u_{xx}$ the second order derivative with respect to $x$. When
$h\rightarrow 0$, instead of forcing the solution of a problem to be
convex, we force it to some mesh-dependent combination of its second
order derivatives to have a sign or another. In addition, this
combination is meaningless since for the convex limit function, it may
have whatever sign. We will say that such a discretization is not
consistent.


%
%

We may repeat here that our non-consistency result is not
contradictory with the convergence results of \cite{Waschmuth_17} (for
instance) since these results rely on a regularization that makes the
algorithm look for strictly convex functions. Looking in the interior
of the set of convex functions, the constraint may not be saturated
and so disappears.

Below, we investigate the consistency of various other discretizations.

\subsection{Use of a weak $P_1$ version}
\label{sec3.2}

We use here a weak $P_1$ definition of convexity which is similar to
the one of \cite{Aguilera_Morin_09}, except that their trial basis and
test basis functions are different. We use the same trial and test
$P_1$ basis. Their method is fully discussed and illustrated in
Section \ref{sec.5}.

\subsubsection{The subharmonicity constraint}

Subharmonicity is an interesting property simpler than convexity since
it is only linear. A weak definition of subharmonicity ($\Delta u \geq
0$) is, for any test function $\phi_i$ in the discrete basis
($\otimes$ is defined in (\ref{def_otimes})):
\begin{equation}
\label{eq14}
\langle  \Delta u_h,\phi_i \rangle = \mbox{Tr }\langle  \mathbf{D}^2 u_h, \phi_i\rangle  = -\mbox{Tr } \langle  \mathbf{\nabla} u_h \otimes \mathbf{\nabla} \phi_i \rangle \geq 0.
\end{equation}
One may then prove the following Proposition:

\begin{proposition}
\label{prop.1}
The weak $P_1$ discretization of the subharmonicity constraint on a mesh like
in Figure \ref{fig3} is consistent:
\begin{equation}
\label{eq15}
\langle  \Delta u_h,\phi_i \rangle = -4u_1+u_2+u_5+u_7+u_4 = (u_{xx}+u_{yy})_{x=x_1,y=y_1}h^2 + O(h^3),
\end{equation}
where $i$ is the node at the center of the cell's group $(x_1,y_1)$.
\end{proposition}

The proof is very easy and left to the reader. Indeed, it is only the
discretization of the Laplacian which is known to be consistent and
even convergent !\\

\begin{figure}[hbtp]
    \centering
    \includegraphics[height=6cm]{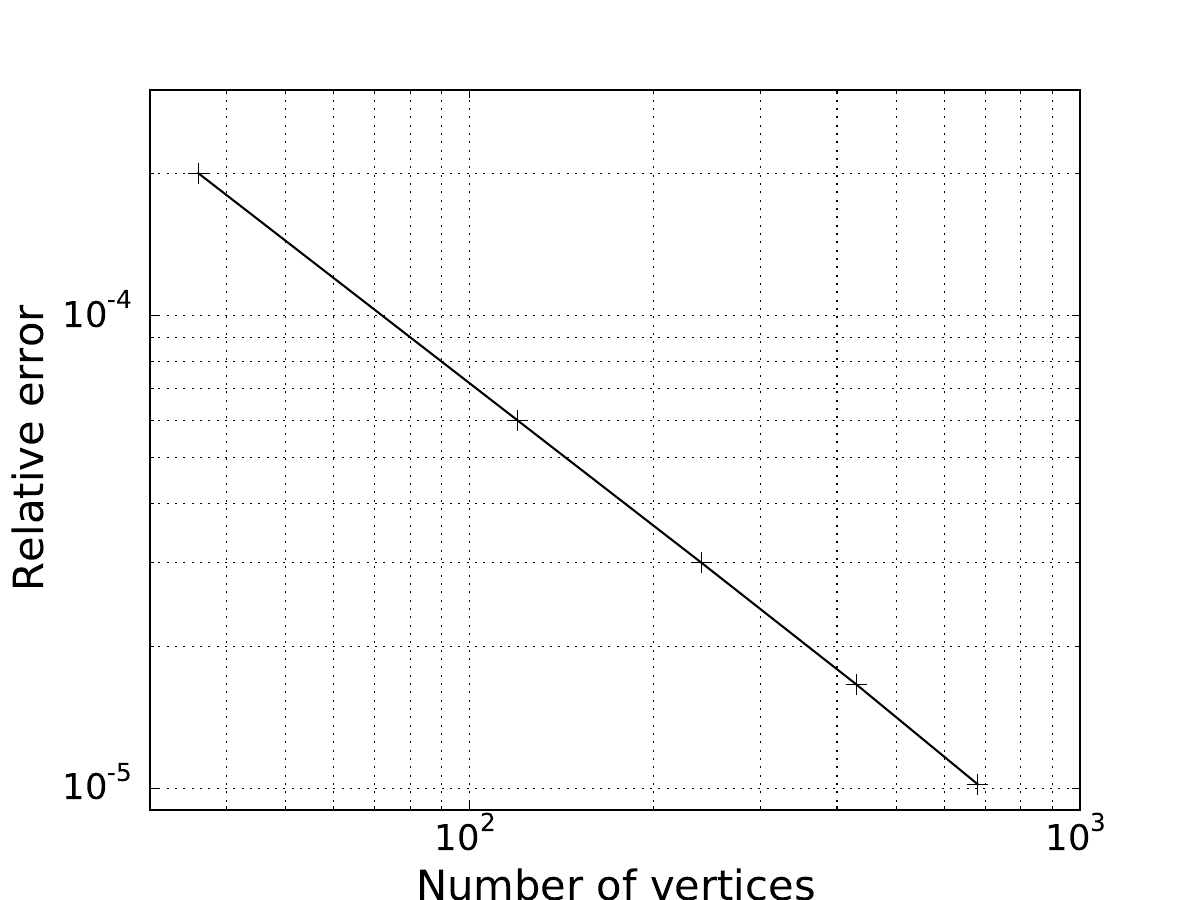}
    \caption{Convergence of $P_1$ FE in case of a subharmonicity constraint.}
  \label{fig3.5}
\end{figure}

In order to test this discretization, we used the Matlab package
\texttt{optim} to minimize the functional $\int_{\Omega} | \nabla u
|^2 /2 + \nabla u_{exact} \cdot \nabla u$, where $f= \Delta
u_{exact}$, over the set of subharmonic functions. The initial
condition is $x(x-1)y(y-1)$. The exact solution of this $H^1_0$
projection is $u_{exact}=x^2y^2-(x^4+y^4)/6$ and its Laplacian
vanishes. The convergence with the procedure \texttt{quadprog} of
quadratic programming can be seen on Figure \ref{fig3.5}. It is quite
satisfactory.


\subsubsection{The convexity constraint}

One may give a weak definition of convexity:

\begin{equation}
\label{eq16}
\mbox{Tr } \langle  \mathbf{D}^2 u_h,\phi_i \rangle  \geq 0 \mbox{ and } \det \langle  \mathbf{D}^2 u_h,\phi_i \rangle  \geq 0,
\end{equation}
for any $\phi_i$, basis function of the FEM at the node $i$. One may
then prove the following Proposition:
\begin{proposition}
\label{prop.2}
The weak $P_1$ discretization of the convexity constraint (\ref{eq16}) on a
mesh as in Figure \ref{fig3} is consistent as can be seen from (\ref{eq15})
and ($\otimes$ is defined in (\ref{def_otimes})):
\begin{equation}
\label{eq17}
\begin{array}{l}
\det \langle  \mathbf{D}^2 u_h,\phi_i \rangle  = \det \langle  \mathbf{\nabla} u_h \otimes \mathbf{\nabla} \phi_i \rangle  = \\
\hspace*{1cm} =(-2u_1+u_2+u_5)(u_7+u_4-2u_1)-(-u_7+u_3-u_4+u_6-u_5-u_2+2u_1)^2 /4 \\
\hspace*{1cm} =(u_{xx}u_{yy}-u_{xy}u_{yx})_{x=x_1,y=y_1}h^4+O(h^5).
\end{array}
\end{equation}
where $i$ is the node at the center of the cell's group (point of
coordinates $(x_1,y_1)$.
\end{proposition}

This discretization (\ref{eq16}) for the (nonlinear) convexity contraint is
consistent. Nevertheless, the numerical treatment of the linear and
nonlinear constraints should raise inaccuracies since they are of very
different orders of magnitude ($O(h^2)$ and $O(h^4)$). A very natural question would be to
investigate its convergence. Such a study is postponed to an other
article.

The consistency of the $P_1$ FEM is tested for various meshes. The
weak $P_1$ Hessian is reported in Table \ref{table1} where one may
claim that two meshes are non-consistent with the determinant but all
are consistent with the trace (laplacian).

\section{Consistency of the $P_2$ FEM}
\label{sec.4}

In the present section, we investigate the consistency of various
discretizations through $P_2$ FEM of two second order derivative
constraints: subharmonicity and convexity.

\begin{figure}[hbtp]
    \centering
    \includegraphics[height=6cm]{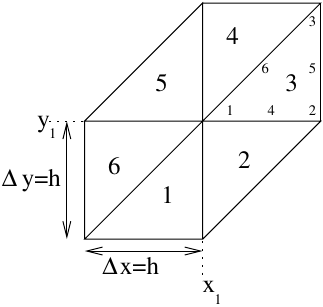}
    \caption{Local shape of the mesh close to $(x_1,y_1)$. Local numbering in triangle 3.}
  \label{fig4}
\end{figure}

First, we interpolate a continuous function $u$ to a $P_2$ function
$u_h= \sum_i u_i \phi_i(x)$ in a domain $\Omega$ meshed with triangles
as in Figure \ref{fig4}. Here, the index $i$ denotes both vertices and
edge midpoints. Then we compute the discretized version of the second
order term constrained to be nonnegative (various versions are
treated) and compute its series expansion. If forcing it to be
nonnegative amounts to forcing the continuous limit function to the
correct constraint, then we claim the discretization is
consistent. Otherwise it is
inconsistent.\\

For the whole section, we assume the mesh is (at least locally)
structured around the point of coordinates $(x_1,y_1)$. The local
numbering of triangles is depicted in Figure \ref{fig4}. The node
$(x_1,y_1)$ is localy numbered 1 in every triangle and the local
numbering of nodes is depicted in triangle $3$ of Figure \ref{fig4}.
%

\subsection{Gradients' jumps for the convexity constraint}
\label{sec4.1}

In a way similar to the $P_1$ case, one may prove for $P_2$ functions
$u_h$:
\begin{align*}
\langle \frac{\partial^2 u_h}{\partial \va \partial \vb}, \varphi \rangle =  \sum_e \int_e (\mathbf{q}_2(s)-\mathbf{q}_1(s)) \cdot \mathbf{n}\,  \varphi(s){\rm d} s \, (\mathbf{n}\cdot \va)(\mathbf{n}\cdot \vb) +\sum_{K} \frac{\partial^2 u_h}{\partial \va \partial \vb}|_{K} \int_K \varphi,
\end{align*}
for any $e$ interior edge of the mesh, $\va$ and $\vb$ are two unit
normal vectors, $\varphi$ is a ${C}_0^{\infty}$ function with compact
support in $\Omega$. By taking $\varphi$ localized along the edge $e$,
one may state that this definition of convexity (with
$\mathcal{C}_0^{\infty}$ test functions) implies the non-negativity of
the gradients' jumps.

In order to test the consistency of such a discretization, we compute
the gradients' jumps accross the edges common to triangles $1$ and
$2$, $2$ and $3$ and last $3$ and $4$. Of course, they are not
constant since they are $P_1$ FE. After an exact computation and a
series expansion (details left to the reader), one may state the following
Proposition.

\begin{proposition}
\label{prop.3}
The discretization of the convexity constraint with the jump of the
gradients between triangles $1$ and $2$, $2$ and $3$, $3$ and $4$ of a $P_2$
function on a mesh like in Figure \ref{fig4} gives terms
\begin{equation}
\label{eq18}
\begin{array}{rcl}
\mbox{Jump(1/2)} & = & (u_{xxy}+u_{xyy})_{x=x_1,y=y_1}(y-y_1+h/2)h/2 + O(h^3); \, y\in [y_1-h,y_1] \\
\mbox{Jump(2/3)} & = & (u_{xxy}+u_{xyy})_{x=x_1,y=y_1}(x-x_1-h/2)h/2 + O(h^3); \, x\in [x_1,x_1+h] \\
\mbox{Jump(3/4)} & = & -(u_{xxy}+u_{xyy})_{x=x_1,y=y_1}(x-x_1-h/2)\sqrt{2}h/2 + O(h^3);\, x\in [x_1,x_1+h].
\end{array}
\end{equation}
Such a discretization is non-consistent.
\end{proposition}

Since for instance between triangles $1$ and $2$, $(y-y_1+h/2)$
changes sign (but remains $O(h)$), one deduces from (\ref{eq18}) that
forcing the non-negativity of the gradients' jump along the whole edge forces
the limit function to satisfy the non-natural equality condition
$u_{xxy}+u_{xyy}=0$ ! 

So the weak $\mathcal{C}_0^{\infty}$ (or strong !) definition of
convexity implies the non-negativity of gradient's jumps which is
non-consistent and so must be rejected. We check below that the weak
$P_2$ definition may (subsection
\ref{sec4.3}) or maynot be consistent (subsection
\ref{sec4.2}).

\subsection{Weak $P_2$ version of the second derivative at a vertex}
\label{sec4.2}

We define the weak $P_2-P_2$ ({\em i.e.} with $P_2$ trial and test functions)
version of the Hessian matrix at vertex $i$ as:
\begin{equation}
\label{eq19}
\langle \mathbf{D}^2 u_h, \phi_i \rangle  = -\langle \mathbf{\nabla}u_h \otimes \mathbf{\nabla} \phi_i \rangle,
\end{equation}
where $\otimes$ is defined in (\ref{def_otimes}). Strictly speaking,
this may not provide a correct weak version of non-negativity since
$\phi_i$ changes sign. Anyway, one may state the following Proposition
which proof is left to the reader:

\begin{proposition}
\label{prop.4}
The discretization of the linear part of the convexity constraint
according to (\ref{eq19}) of a $P_2$ function on a mesh as in Figure
\ref{fig4} ($h=\Delta x =\Delta y$) and a $P_2$ test function centered
at the vertex $(x_1,y_1)$ gives:
\begin{equation}
\label{eq20}
\mbox{Tr }(\langle \mathbf{D}^2 u_h, \phi_i\rangle ) = -(u_{xxxx}+u_{yyyy})_{x=x_1,y=y_1} * h^4 /48 + O(h^5).
\end{equation}
Such a discretization of an inequality constraint is non-consistent.
\end{proposition}

It appears that such a discretization is not even consistent for the
linear part of the constraint. More precisely, while one could believe
one forces the solution to be subharmonic, indeed, one forces it to be
such that $u_{xxxx}+u_{yyyy} \leq 0$ ! The full nonlinear convexity
constraint on the same vertices may only work worse.

One must notice that the fourth order of the expansion in (\ref{eq20})
is meaningful. Indeed, on the one hand the three midpoints quadrature
is exact for $P_2$ functions in a triangle. On the other hand the
basis function for a vertex vanishes on these midpoints. As a
result, the order two term is identically zero. So the first
non-zero term is the fourth one and the constraint is of fourth order
too.

\subsection{Weak $P_2$ version of the second derivative at an edge midpoint}

\label{sec4.3}

We define the weak $P_2$ version of the second derivative at an edge
midpoint (denoted by index $j$) in a way similar to
(\ref{eq19}). Indeed, the function $\phi_i$ in (\ref{eq19}) is
replaced by the $P_2$ basis function $\phi_j$ associated to an edge midpoint
indexed by $j$.

As in subsection \ref{sec3.2}, the basis functions (of the edge midpoints)
are non-negative. So it is {\em a priori} an admissible weak formulation
of semi-definite positiveness.

\subsubsection{The subharmonicity constraint}

Let us assume we discretize the constraint $\Delta u \geq 0$ by

\begin{equation}
\label{eq20.5}
\int_{\Omega} \nabla u_h \cdot   \nabla \phi_j \leq 0,
\end{equation}
for all $j$ index of an edge midpoint interior to $\Omega$. One may
then state a Proposition (which proof is left to the reader):

\begin{proposition}
\label{prop.5}
The discretization of the linear part of the convexity constraint
according to (\ref{eq20.5}) on a mesh as in Figure \ref{fig4}
($h=\Delta x =\Delta y$) gives the same series expansion, whether the
edge is vertical, horizontal or diagonal:
\begin{equation}
\label{eq21}
\mbox{Tr } (\langle \mathbf{D}^2 u_h, \phi_j\rangle ) = (\Delta u)_{x=x_1,y=y_1}h^2/3+O(h^3),
\end{equation}
for any $j$ index of an edge midpoint interior to $\Omega$. Such a
discretization of an inequality constraint is consistent.
\end{proposition}

Various questions remain. If we discretize the subharmonicity (or the
convexity) constraint only at edge midpoints, is it enough constraints
or not ? More generaly, the amount of constraints compared to the
amount of degrees of freedom (dof) should be discussed. The relevance
of this question seems to be the conclusion of most recent articles
like \cite{Mirebeau_14,Merigot_Oudet_14,Aguilera_Morin_09}.



\subsubsection{The convexity constraint}

Like in the subharmonic case, we take a weak $P_2$ version of the
continuous nonlinear constraint det $D^2 u \geq 0$ with (nonnegative)
test functions associated to every edge midpoint $j$ in the
interior. One may then easily prove the following Proposition:

\begin{proposition}
\label{prop.6}
The discretization of the nonlinear part of the convexity constraint
on a mesh as in Figure \ref{fig4} ($h=\Delta x
=\Delta y$) gives:
\begin{equation}
\label{eq22}
\det (\langle \mathbf{D}^2 u_h, \phi_j\rangle ) =(u_{xx}u_{yy}-u_{xy}u_{yx})_{x=x_1,y=y_1} h^4/9 + O(h^5),
\end{equation}
for any $j$ index of an edge midpoint interior to $\Omega$. Such a
discretization of an inequality constraint is consistent.
\end{proposition}

So the convergence order of the linear constraint is two while the
nonlinear one's is four. Such a discrepancy between the orders of
convergence of those two constraints should be managed in numerical
simulation using both constraints. At least this discretization is
consistent.

\section{Discussion on the Aguilera and Morin's articles}
\label{sec.5}

Let us notice that \cite{PCLM} seems not to have been known of the
authors of \cite{Aguilera_Morin_08,Aguilera_Morin_09}.

Roughly speaking, the first article \cite{Aguilera_Morin_08} proves
convergence of the FD discretization of problems like ours. The second
\cite{Aguilera_Morin_09} proves convergence of the FE discretization
of the same problems but the authors imagine the clever trick of using
different basis for the trial and the test functions.

\subsection{Finite Differences}

In \cite{Aguilera_Morin_08}, Aguilera and Morin deal with FD and
``prove convergence under very general conditions, even when the
continuous solution is not smooth''. Their proof relies on
approximation theory and on a weak definition of convexity, namely the
FD approximation of the Hessian matrix is forced to be positive. They
stress that their FD-convexity is not equivalent to a continuous
convexity. As a consequence, their approximation is {\em not
  conformal}. So there is no contradiction with \cite{PCLM}. In
addition, one must notice that the number of constraints is roughly
twice the number of interior points. It is large, but grows only
linearly with the degrees of freedom. They also provide numerical
experiments.

\begin{figure}[hbtp]
    \centering
   \includegraphics[height=5cm]{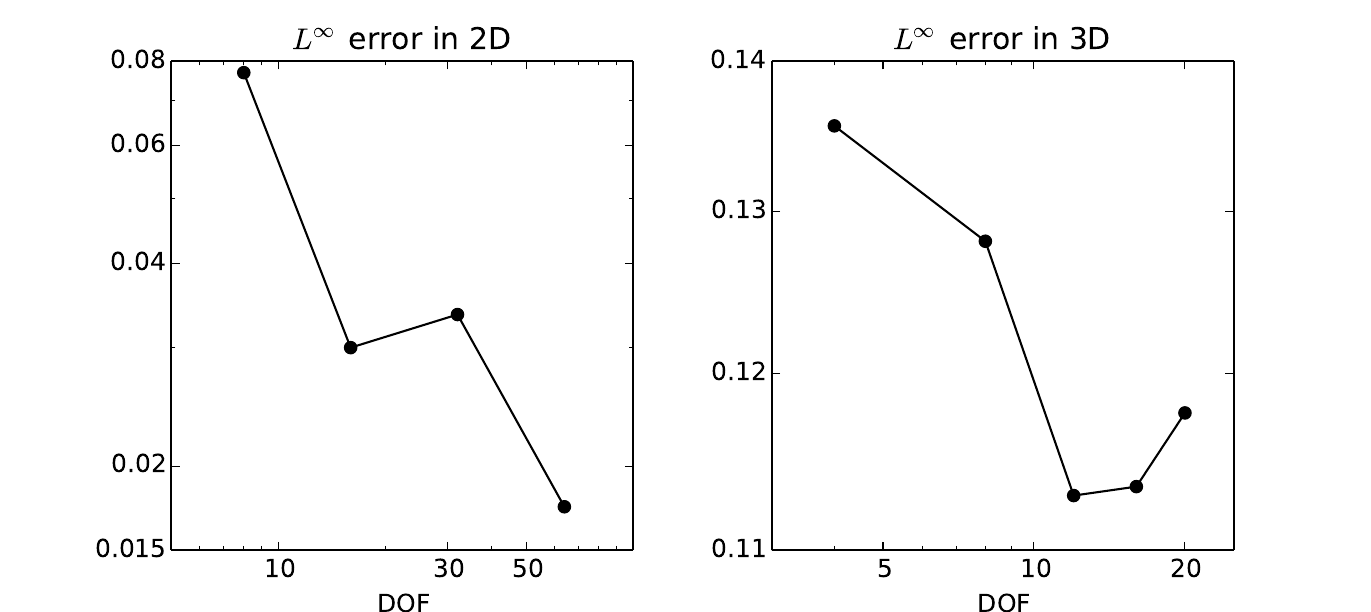}
    \caption{The monopolist numerical solution error versus DOF.}
  \label{fig_AM_1}
\end{figure}

Unfortunately the numerical experiments are not very
convincing. Concerning the monopolist problem in 2D, they claim that
``the error in the $L^{\infty}$ norm is smaller or approximately equal
to $h$'' (p. 27). Yet, their Table 1 provides errors that we draw in a
loglog scale in the left part of Figure \ref{fig_AM_1}. The order of
convergence is not clear.

For the monopolist problem in 3D, the authors notice that ``the
$L^{\infty}$ error is not converging to zero with order $O(h)$''
(p. 29). Since they only prove convergence and do not forecast any
order, there is no contradiction. But when their Table 2 is drawn in a
loglog scale (see Figure \ref{fig_AM_1} right), convergence seems not
to be reached yet. Moreover the execution time grows faster than
polynomialy both for 2D and 3D monopolist.

As reminded by the authors, FD discretization gives the very same
matrices as the $P_1$ FE (both for trial and test functions), on a
mesh like our mesh 1 on Figure \ref{fig1} (except boundary conditions
irrelevant here).  In this sense, one may claim that their FD
discretization is somehow equivalent to a $P_1$ (for trial) - $P_1$
(for test) discretization which was proved to be
non-consistent. Indeed, the next article \cite{Aguilera_Morin_09}
gives one example (Example 3.7 p. 3150) of $P_1$ discretization and a
weak definition of convexity (with $P_1$ trial and test functions, so
alike FD discretization), where `` $(u_h)$ does not converge to $u$ as
$h \rightarrow 0$'' (p. 3151 and Figure 3.1).

Discussing the way these results on FD are coherent with the results
on FE is postponed to a forthcoming article.



\subsection{Finite Elements}

In \cite{Aguilera_Morin_09}, Aguilera and Morin prove convergence of
the discretization of the full problem: not only the approximation of
the Hessian, but also of the functional together.

Especially, they define $u_h$, interpolation of $u$ on a given mesh
with respect to the trial functions basis $\{\phi_r^h\}_r$ as
$u_h(x)=\sum_{r \in I_{trial}^h} u_{hr} \phi_r^h(x)$. They also use a
(possibly) different FE basis for the test functions
$\{\varphi_s^h\}_{s \in I_{test}^h}$. Then, the discrete Hessian matrix is
defined weakly by $H_s^hu_h=-\left( \langle \partial_i u_h, \partial_j
\varphi_s^h \rangle \right)_{i,j}$ and FE convexity is defined as:
\[
H_s^hu_h=-\left( \sum_{r \in I_{trial}^h} u_{hr} \langle \partial_i \phi_r^h , \partial_j \varphi_s^h \rangle  \right)_{i,j}  \succeq 0,
\]
for all $\varphi_s^h$ in the test functions basis, where they denote
$\mathbf{A} \succeq 0$ if $\mathbf{A}$ is positive
semi-definite. Notice that the matrix $(\langle \partial_i \phi_r^h
, \partial_j \varphi_s^h \rangle)_{r,s}$ is a rectangle if the trial
and test basis do not have the same number of functions !

Then their main result states convergence of the discretization once
one assumes at least (p. 3147):
\begin{itemize}
\item[C.2] There exists a linear operator $\mathcal{I}^h$ with
  values in the discretization space $V_h$ (the interpolant), an
  integer $m \geq 2$, and a constant $C$ independent of $u$ and $h$
  such that
\[
\parallel u-\mathcal{I}^hu \parallel_{H^1(\Omega)} \leq C
h^{m-1} \parallel u\parallel_{H^m(\Omega)},
\]
\item[C.3] The basis test-functions are such that
\[
\varphi_s^h(x) \ge 0 \text{ for all } x\in \Omega, \text{ all } h>0, \text{and all }s.
\]
\end{itemize}
The condition C.3 merely states that the test functions, that are
supposed to approximate positive $\mathcal{C}_0^{\infty}$ functions,
remain nonnegative. So the basis functions are at most $P_1$ or a
subset of $P_k$.\\

The condition C.2 is very classical and states that the trial space is
$P_{m-1}$. In addition, the proof requires $m>2$. In other words,
the trial functions must at least be $P_2$.

This is seen by the authors who say they ``need to elaborate'' on the
condition $m>2$. On their p. 3150, they notice the FE discretized
Hessian ($P_1-P_1$ and so $m=2$) expands into the FD scheme (for one
specific mesh !) at first order and so they refer to their FD article
to conclude that in the case of some specific meshes, their Theorem
``also holds for $m=2$'' (or $P_1-P_1$ FE). \\

Yet, the very next example they provide is $P_1-P_1$ and they ``report
numerical evidence supporting the necessity of assuming $m>2$ in
[their] Theorem 3.6'' (p. 3150). They stress that ``since $u$ is
convex, the projection $u_h$ should converge to $u$ as $h \rightarrow
0$, but this is not the case in this example'' (p. 3150) and ``$(u_h)$
does not converge to $u$ as $h \rightarrow 0$'' (p. 3151). They
summarize this result by saying ``although there is some sort of
super-convergence for some meshes, for general meshes [...] FE-convex
piecewise {\em linear} [m=2!] function may not suffice''. Discussing
how this numerical example articulates with \cite{Aguilera_Morin_08}
is postponed to a forthcoming article.\\

All this is even complexified when the authors cope with numerical
experiments of $P_2-P_2$ FE discretization. As they notice, the FE
test functions in the basis are not all nonnegative. So they
``considered the usual piecewise {\em linear} nodal basis for the
vertices'' and ``the usual quadratic bubbles'' for the edges'
midpoints (p. 3152).

\begin{figure}[hbtp]
    \centering
    \includegraphics[height=5cm]{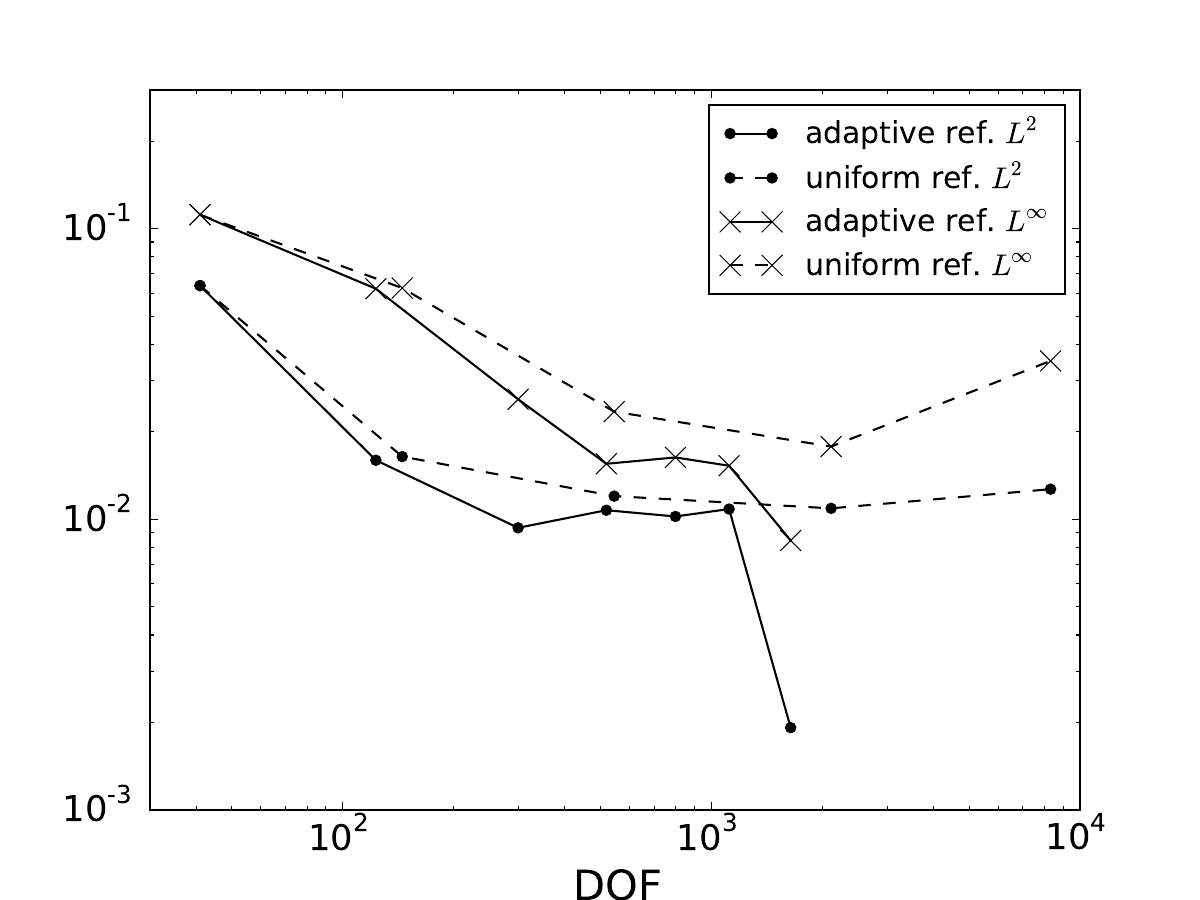}
    \caption{Error according to the refinement either adaptive (solid lines) or uniform (dashed lines) and the norm ($L^2$ or $L^{\infty}$) versus Degrees of Freedom (DOF).}
  \label{fig_AM_2}
\end{figure}

Their table of convergence (Table 5.1) is transformed into loglog
graphic in Figure \ref{fig_AM_2}. The last simulation with {\em
  adaptive} refinement gives a good point, but the uniform-refinement
simulations give non-decreasing errors. Again convergence is not
obvious. So as to make this convergence obvious, we illustrate their
method in Section \ref{sec5.4}.




\subsection{Are \cite{Aguilera_Morin_08,Aguilera_Morin_09} contradictory with \cite{PCLM} ?}
\label{AM_moi}

The answer is no. But it deserves to be explained.

Both \cite{Aguilera_Morin_09} and \cite{PCLM} use a dual definition of
convexity and have conclusions that could seem to be
contradictory. The only difference is that \cite{Aguilera_Morin_09}
uses different basis for trial and test functions:
\[
-\langle \nabla u_h \otimes \nabla \varphi_s^h \rangle  \succeq 0,
\]
while \cite{PCLM} uses the same basis for trial and test functions,
both in $P_1$:
\[
-\langle \nabla u_h \otimes \nabla \phi_h \rangle \succeq 0.
\]
Both discretizations are weak in a sense, but the one of
\cite{Aguilera_Morin_09} has less test functions than trial functions
and so less contraints than the number of Degrees Of Freedom
(DOF). Indeed $m >2$ means that the trial basis is rather large (at
least $P_2$) and condition C.3 that the test basis is at most
$P_1$. The matrix of constraints is rectangle and there are less
constraints than the unknowns whereas the discretization of
\cite{PCLM} uses the same basis for trial and test and so involves
square matrices.

Then, it is not surprising that the problem discretized with the weak
$P_1-P_1$ definition in \cite{PCLM} is overconstrained : there are as
many constraints as DOF ! Indeed it appears also from the literature
that the amount of constraints must be not too large.  The article
\cite{Mirebeau_14} defines a stencil of constraints that is coarser
than the grid and \cite{Merigot_Oudet_14} defines a subgrid of points
on which the constraint is applied. Then it seems coherent that the FE
proof of \cite{Aguilera_Morin_09} requires {\em more} DOF than the
number of constraints. The methods of \cite{Waschmuth_17} and similar
works look for functions in the interior of the convex
functions set where the boundary  are irrelevant.\\

Our Proposition \ref{prop.4} proved in subsection \ref{sec3.2} states that
if we discretize $u_h$ with $P_2$ FE and use a $P_2$ set of test
functions (but then some test functions change sign), then
the discretization of the linear part of the constraint is
non-consistent. It proves that the proof of convergence of FEM
by \cite{Aguilera_Morin_09}, which fails if the test functions change
sign (as in $P_2$ and higher), maynot be improved on that point. So we do believe their theorem is optimal.

If we use the consistency test for some discretizations, we can
compute the first non-zero term in the expansion of
\[
-\langle \nabla u_h \otimes \nabla \varphi_s^h \rangle,
\]
for $P_1-P_1$ discretizations and report these computations in Table
\ref{table1} (lower part) on various meshes depicted in Table
\ref{table1} (upper part). It appears that some discretizations are
not consistent with convexity (two on the left) while others are (two
on the right).

\subsection{Does the $P2-P_1$ FEM converge ?}
\label{sec5.4}

The convergence is proved in \cite{Aguilera_Morin_09} for $P_k$
($k\geq 2$) trial functions and $P_1$ test functions with additionnal
assumptions. In order to ensure the numerical convergence, we wrote a
Python code to mesh a square $\Omega=[1,2]^2$, discretize the
functional (\ref{eqm3}), and discretize the nonlinear constraints of
inequality (convexity). Then we use the \verb+minimize+ function in
\verb+scipy+ to minimize the functional. The execution time are
meaningless since they rely on the level of precision we ask and we
only want to illustrate convergence.

In a very first step, we look for the solution to
(\ref{eqm3},\ref{eqm2},\ref{eq5}), which exact solution is given in
(\ref{eq5.5}). It is not surprising that, {\em whatever the mesh}, we
can get a relative error of $10^{-12}$. It only justifies that a
quadratic function may be approximated by a $P_2$ numerical code without error.

In order to justify convergence, one needs a more complex convex
function to find.

\subsubsection{Convergence to an exact solution}

We choose
\[
u_{exact}(x) = \exp((x_1+x_2)/2)-e,
\]
and take the functional
\[
\int_{\Omega} \mid \nabla u \mid^2 + \Delta u_{exact} u \; {\rm d}x-\int_{\partial \Omega} \DP{u_{exact}}{n} u ,
\]
which is also, up to a constant, $\int_{\Omega} \mid \nabla
(u-u_{exact}) \mid^2$. The function is constrained to be convex and be
such that $u(1,1)=0$, which makes sense since the function is convex
and so regular enough for evaluating at $(1,1)$. The solution is
convex, but we are going to prove that its $P_2$ interpolate is {\em
  not} convex.

\begin{figure}[hbtp]
    \centering
    \includegraphics[height=5cm]{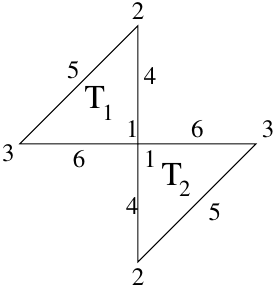}
    \caption{Local numbering.}
  \label{fig6}
\end{figure}

The geometrical argument is that we $P_2$-interpolate a function of
$x_1+x_2$ on a mesh like mesh 1 in Figure \ref{fig1}. Let us take the
$P_2$ interpolate of $u_{exact}$ on the two triangles depicted in
Figure \ref{fig6}, close to the point $(x_1,y_1)$ at the center as before. Let
us prove now that the $P_2$ interpolate $u_h^{P_2}=\Pi_{P_2}u \in P_2$
is strictly convex inside $T_2$. Along the segment between the point
numbered 1 in $T_2$ and the point 5 in $T_2$, one has:
\[
\begin{array}{rl}
u_h^{P_2}(x_1+\alpha h,y_1-\alpha h)=&2u_1^{T_2}(-2\alpha+1)(-2\alpha+1/2)+2u_2^{T_2}\alpha(-\alpha +1/2)\\
 & +2u_3^{T_2}\alpha (\alpha-1/2)+4 u_4^{T_2} \alpha(-2\alpha+1)\\
 & +4u_5^{T_2}\alpha^2+4u_6^{T_2}\alpha(-2\alpha+1)\\
=& e^{(x_1+y_1)/2}\left( \alpha^2(12+2e^{-h/2}+2e^{h/2}-8e^{-h/4}-8e^{h/4})\right.\\
 & \left.+\alpha(-6-e^{-h/2}-e^{h/2}+4e^{-h/4}+4e^{h/4})+1 \right),
\end{array}
\]
where $u_i^{T_2}$ is the value of the interpolate at the DOF locally
numbered $i$ in triangle $T_2$. These values are the same as those of
$u_{exact}$ by definition of the $P_2$ interpolate and enable to write
the last equality. It is then easy to check that this function is
strictly convex for $\alpha \in (0,1/2)$. The same computations in the
triangle $T_1$ gives for $\alpha \in (0,1/2)$
\[
\begin{array}{rl}
u_h^{P_2}(x_1-\alpha h,y_1+\alpha h)=&e^{(x_1+y_1)/2}\left( \alpha^2(12+2e^{-h/2}+2e^{h/2}-8e^{-h/4}-8e^{h/4})\right.\\
 & \left.+\alpha(-6-e^{-h/2}-e^{h/2}+4e^{-h/4}+4e^{h/4})+1 \right).
\end{array}
\]
Since both $u_h^{P_2, T_2}$ and $u_h^{P_2, T_1}$ are {\em strictly} convex
in their domains and have the {\em same} value on three aligned points, the
function $u_h^{P_2}$ maynot be convex at least on the segment between
the points numbered 5 in $T_2$ and the point numbered 5 in $T_1$ (see
Figure \ref{fig6}).

\begin{figure}[hbtp]
    \centering
    \includegraphics[height=5cm]{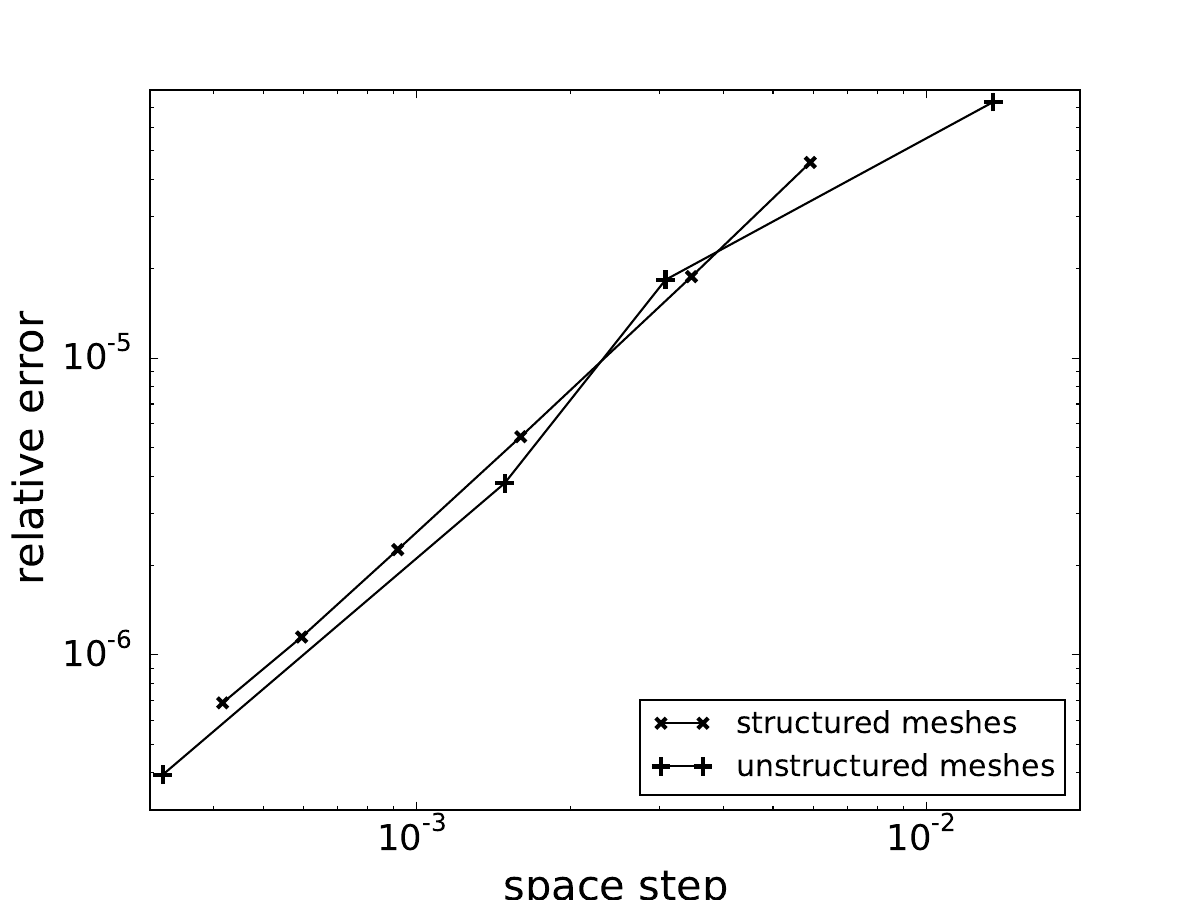}
    \caption{Convergence of the FEM $P_2-P_1$ method.}
  \label{fig7}
\end{figure}

Despite the non-convexity of the $P_2$ interpolate of the exact
solution (which is convex), the convergence of the FEM method proposed
and proved in \cite{Aguilera_Morin_09} is clearly illustrated in
Figure \ref{fig7} where structured and unstructured meshes are used.


The order of convergence seems to be the same on structured and
unstructured meshes. It is roughly 1.6.

\subsubsection{Can this method be used ?}

In order to use this method justified theoretically (in
\cite{Aguilera_Morin_09}) and numerically (above), we look for a
solution to the problem of the monopolist in $\Omega=[1,2]^2$:
\[
\displaystyle \min_{\begin{array}{c}u \; \mbox{convex}\\ u\geq 0\end{array}}\int_{\Omega} \mid \nabla u \mid^2-x \cdot \nabla u + u \, {\rm d}x.
\]
There exists no explicit solution to this problem but some properties
of the solution are known (see \cite{RC}).

\begin{figure}[hbtp]
  \begin{minipage}[c]{.46\linewidth}
    \centering
    \includegraphics[height=8cm]{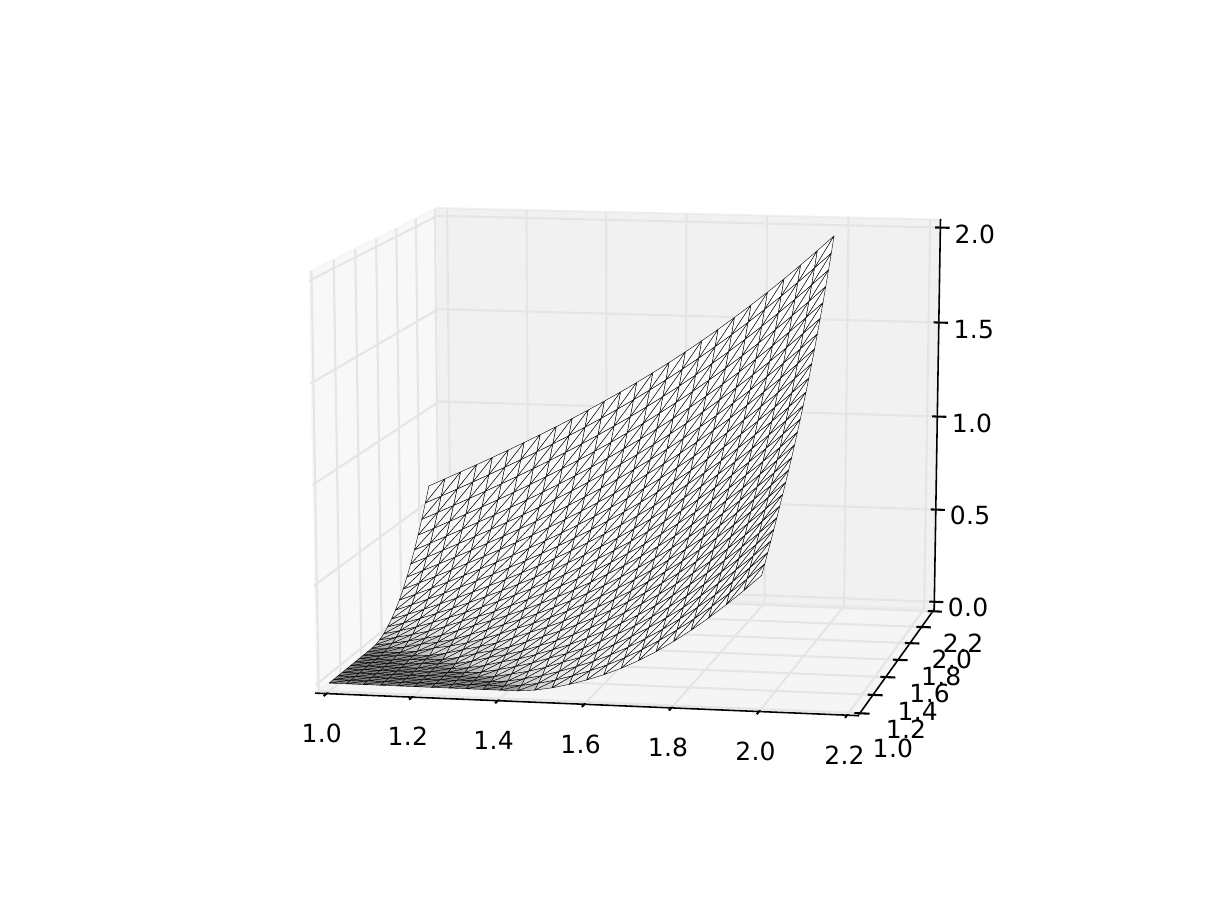}
  \end{minipage}
   \hfill%
   \begin{minipage}[c]{.46\linewidth}
     \centering
     \includegraphics[height=6cm]{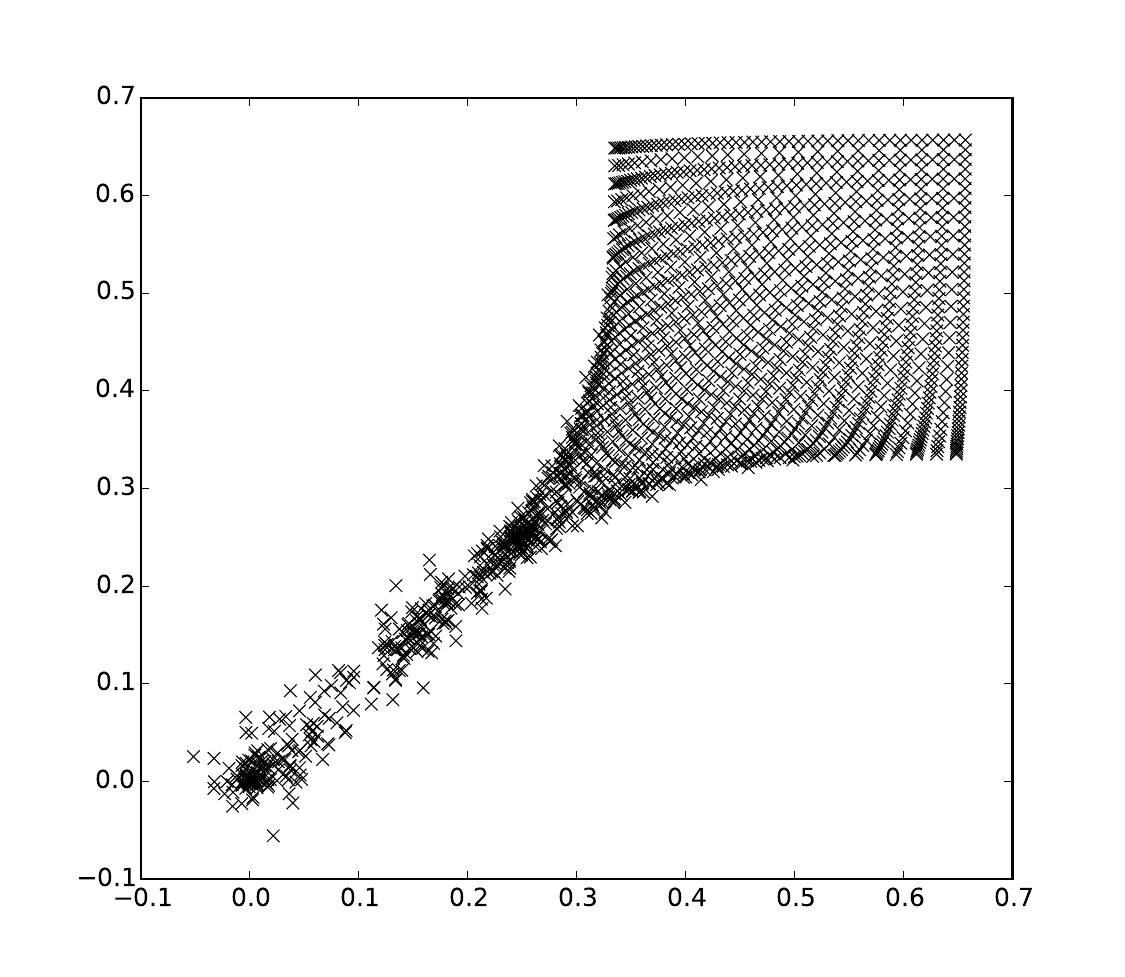}
   \end{minipage}
  \caption{The monopolist's graph and its gradient ($P_2-P_2$) on a $25\times 25$ regular mesh.}
  \label{fig8}
\end{figure}

The computed $P_2$ function on a structured mesh $25 \times 25$, its
gradient and its second eigenvalue look similar to the ones computed
in \cite{Mirebeau_14}. The graph of the full $P_2$ function would be
difficult to understand because of the number of DOF. That is why we
only use the $P_1$ information from $u_h$ in Figure \ref{fig8} (left).

So as to display the gradient, one might compute the gradient
$\int_{\Omega} \nabla u_h \phi_i$ for $i$ associated to an interior
point, with a $P_2$ function $u_h$, and $\phi_i \in P_2$. But since
the quadrature uses the three mid-points, it is right only up to $P_2$
functions, and so the integral being over a $P_3$ function is not
exact. Yet, since the $P_2-P_1$ discretization, which is exact, gives
the same graphical result with less points, we consider the Figure
with $P_2$ test functions to be acceptable. The result may be seen in
Figure \ref{fig8} (right).

\begin{figure}[hbtp]
  \begin{minipage}[c]{.46\linewidth}
    \centering
    \includegraphics[height=7cm]{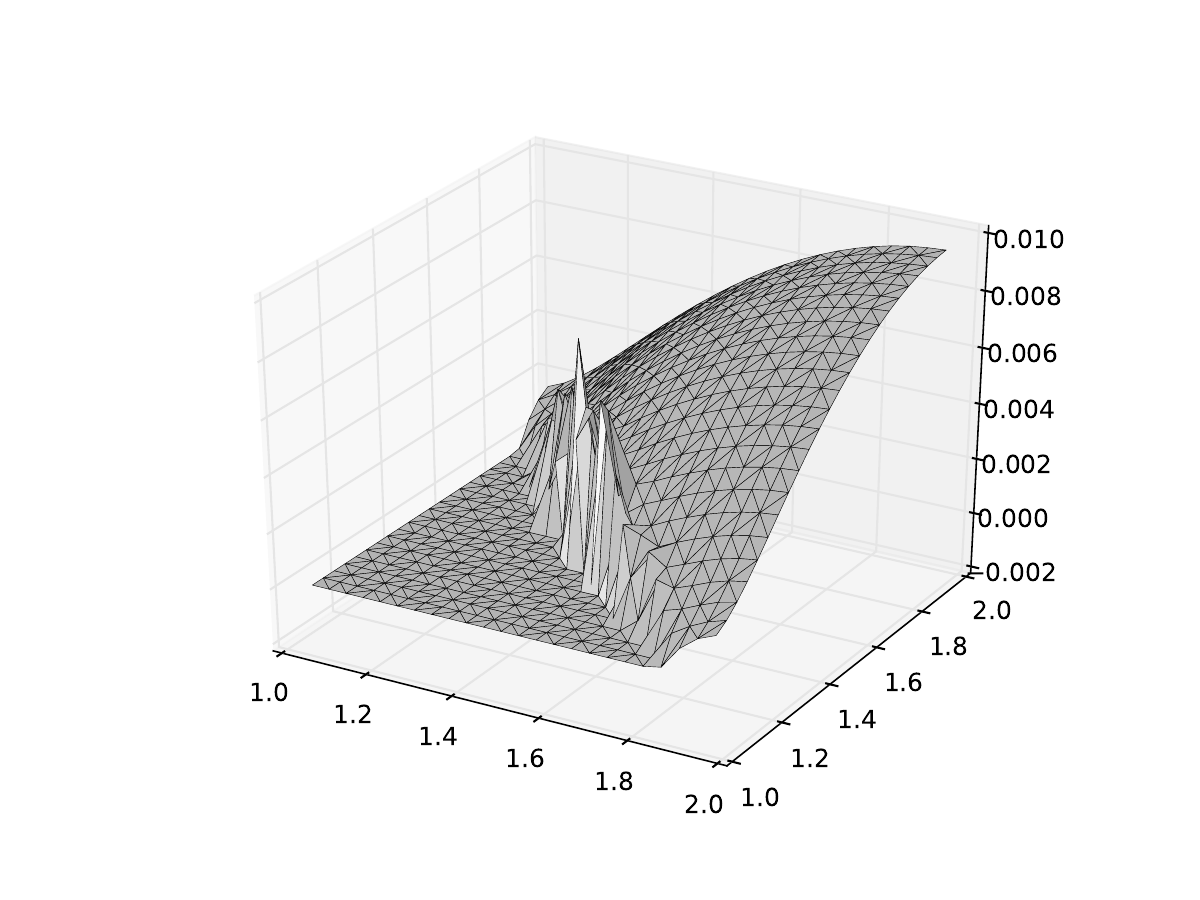}
  \end{minipage}
  \hfill%
  \begin{minipage}[c]{.48\linewidth}
    \centering
    \includegraphics[height=6cm]{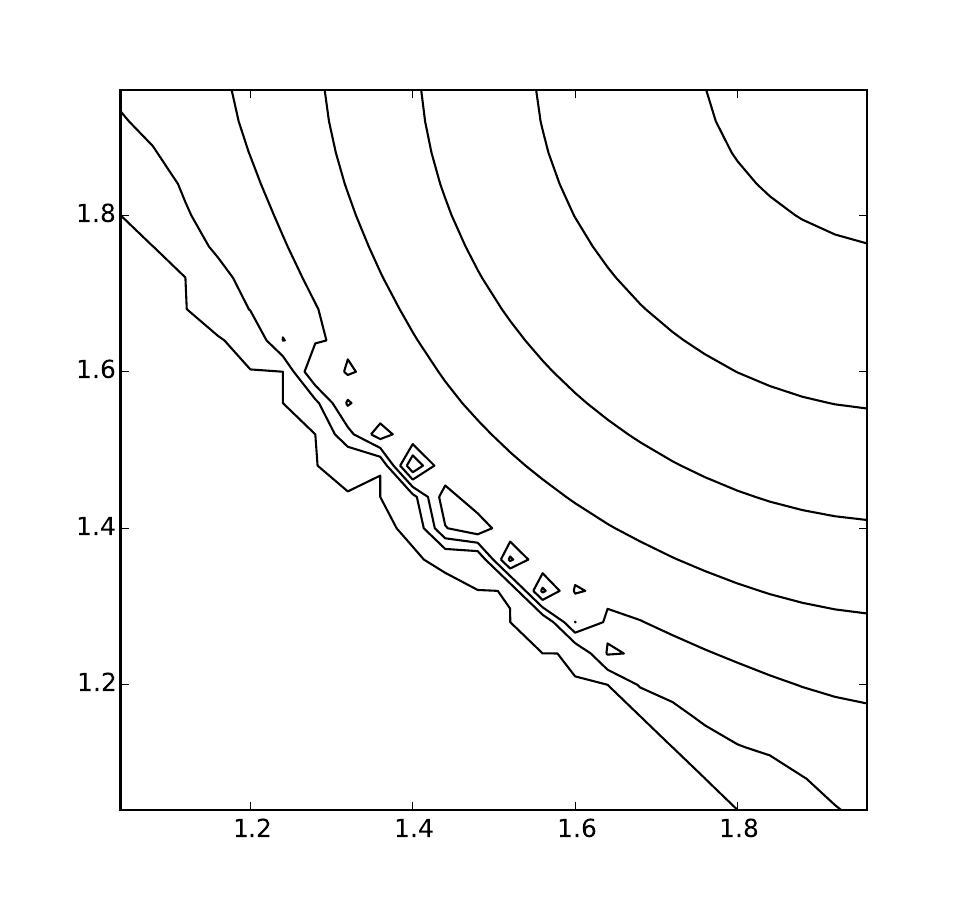}
  \end{minipage}
  \caption{The eigenvalues.}
  \label{fig9}
\end{figure}

In Figure \ref{fig9}, one may see the second eigenvalue of the
Hessian, computed as a $P_2$ function tested on $P_1$ functions and
its level sets on the right. The first eigenvalue is not significant
because of numerical artifacts. They are computed as the two roots per
interior vertex of:
\[
-\left(
\begin{array}{cc}
\int_{\Omega} \partial u_h / \partial x \; \partial \phi_i / \partial x & \partial u_h / \partial x \; \partial \phi_i / \partial y \\ 
\int_{\Omega} \partial u_h / \partial y \; \partial \phi_i / \partial x & \partial u_h / \partial y \; \partial \phi_i / \partial y
\end{array}
\right).
\]


\section{CONCLUSION}
\label{sec.6}

In this article we prove (Theorem \ref{th6}) that the $P_1$
discretization of a function that satisfies a strong definition of
convexity (with $\mathcal{C}_0^{\infty}$ test functions), which is
equivalent to the gradients' jumps positivity (for $P_1$ FE), and so is
conformal, leads to an additional constraint on the limit
function. While \cite{PCLM} dealt mainly with regular meshes, our result
extends explicitely to any mesh. The error of such a discretization of the
constraint does not vanish with the space step. We justify that it is
localized where the additional constraint on the limit function is not
satisfied. The condition for existence of such a counter-example
requires information both from the mesh and its refinement but is
general.

In addition, among the $P_1$ discretization of $u$ with a strong (section
\ref{sec3.1}) or weak (section \ref{sec3.2}) definition
of convexity, not all are consistent. The definition of consistency is
very similar to the one of partial differential equations and it is
used above to discriminate likely discretizations and unlikely
ones. But this does not guarantee good numerical results even when the
discretization is consistent. We also discuss some structured meshes,
the $P_1$ discretization of $u$ and the use of $P_1$ test functions
for convexity (say $P_1-P_1$ convexity). we prove they may or maynot be consistent
(Table \ref{table1}).

We also test the gradients' jumps of $P_2$ functions (strong
discretization with $\mathcal{C}_0^{\infty}$ test functions in
subsection \ref{sec4.1}) and various weak $P_2$ discretizations for
$u_h$ and for some $P_2$ test functions (subsections \ref{sec4.2} and
\ref{sec4.3}). Some of them are consistent and some are not.

We also discuss the literature, but let a deeper study of how
\cite{Aguilera_Morin_08,Aguilera_Morin_09} interact with \cite{PCLM}
to a forthcoming article. In a last subsection, we numerically
illustrate the convergence of \cite{Aguilera_Morin_09}'s FEM where the
trial and test basis are different. Even if other methods exist, such as the recent methods using an evolution PDE converging to the convexified function \cite{Oberman_08,Carlier_Galichon_11}, this
FEM is maybe one of the simplest. Comparisons of all these methods would
be a thrilling challenge.

\section*{Acknowledgments}
The author thanks G. Carlier for drawing his attention to the problem
of subharmonicity approximation.


\begin{thebibliography}{99}
\bibitem{PCLM} P. Chon\'e and H.V.J. Le Meur (2001). Non-convergence result for conformal approximation of variational problems subject to a convexity constraint. {\em Numer. Funct. Anal. Optim.} 22 no. 5-6:529--547.
\bibitem{Newton} I. Newton (1686). {\em Philosophiae Naturalis Principia Mathematica}.
\bibitem{Goldstine} H.H. Goldstine (1980). {\em A history of the calculus of variations from the 17th to the 19th century} Springer-Verlag, Heidelberg.
\bibitem{buttazo} G. Buttazzo, V. Ferone and B. Kawohl (1995). Minimum problems over sets of concave functions and related questions. {\em Math. Nachr.} 173:71-89.
\bibitem{LR_Pelletier_01} T. Lachand-Robert and M.A. Peletier (2001). An Example of Non-convex Minimization and an Application to Newton's Problem of the Body of Least Resistance, {\em Ann. Inst. H. Poincar\'e Anal. Non Lin\'eaire} 18, no. 2:179--198.
\bibitem{PL2} P.L. Lions (1998). Identification du c\^one dual des fonctions convexes et applications.  (French. English, French summary)  [Identification of the dual cone of convex functions and applications] {\em C. R. Acad. Sci. Paris S\'er. I Math.} 326 no. 12, 1385--1390.
\bibitem{LR_Oudet_05} T. Lachand-Robert and \'E. Oudet (2005). Minimizing within convex bodies using a convex hull method, {\em SIAM Journal on Control and Optimization} Vol. 16 Number 2:368-379.
\bibitem{Waschumth_14} G. Wachsmuth (2014). The numerical solution of Newton’s problem of least resistance, {\em Math. Program.} Ser. A 147: 331--350. 
\bibitem{Carlier_04} G. Carlier (2004). On a theorem of Alexandrov. {\em J. Nonlinear Convex Anal.} 5, no. 1:49--58.
\bibitem{Carlier_Comte_Peyre_09} G. Carlier, M. Comte and G. Peyr{\'e} (2009). Approximation of maximal {C}heeger sets by projection, {\em M2AN Math. Model. Numer. Anal.} 43(1):139-150.
\bibitem{RC} J.-C. Rochet and P. Chon\'e (1998). Ironing, sweeping and multidimensionnal screening, {\em Econometrica} {\bf 66}:783-826.
\bibitem{Carlier_Lachand-Robert_01} G. Carlier and T. Lachand-Robert (2001). Regularity of solutions for some variational problems subject to a convexity constraint. {\em Comm. Pure Appl. Math.} 54 no. 5:583--594.
\bibitem{Leung_Renka_99} N.K. Leung and R.J. Renka (1999). $C\sp 1$ convexity-preserving interpolation of scattered data. {\em SIAM J. Sci. Comput.} 20 , no. 5:1732--1752.
\bibitem{Carlier_LR_Maury_99} G. Carlier, T. Lachand-Robert and B. Maury (1999). $H\sp 1$-projection into the set of convex functions: a saddle-point formulation.  {\em CEMRACS 1999 (Orsay)},  277--289 (electronic), {\em ESAIM Proc., 10, Soc. Math. Appl. Indust.}, Paris.
\bibitem{Maury_03} B. Maury (2003). Version continue de l'algorithme d'Uzawa. (French) [Continuous version of the Uzawa algorithm]  {\em C. R. Math. Acad. Sci.} Paris 337 , no. 1:31--36.
\bibitem{Carlier_LR_Maury_01} G. Carlier, T. Lachand-Robert and B. Maury (2001). A numerical approach to variational problems subject to convexity constraint, {\em Numer. Math.} 88, no. 2:299--318.
\bibitem{Aguilera_Morin_08} N.E. Aguilera and P. Morin (2008). Approximating optimization problems over convex functions, {\em Numer. Math.} 111(1):1-34
\bibitem{Aguilera_Morin_09} N.E. Aguilera and P. Morin (2009). On convex functions and the finite element method, {\em SIAM J. Numer. Anal.} 47(4):3139-3157
\bibitem{Ekeland_Moreno_10} I. Ekeland and S. Moreno-Bromberg (2010). An algorithm for computing solutions of variational problems with global convexity constraints, {\em Numer. Math.} 115(1):45-69.
\bibitem{Waschmuth_17} G. Wachsmuth (2017). Conforming Approximation of Convex Functions with the Finite Element Method, submited.
\bibitem{Merigot_Oudet_14} Q. M\'erigot and \'E. Oudet (2014). Handling convexity-like constraints in variational problems, {\em SIAM J. Numer. Anal.}, 52 (5):2466--2487.
\bibitem{Mirebeau_14} J.M. Mirebeau (2015). Adaptive, Anisotropic and Hierarchical cones of Discrete Convex functions, {\em Numer. Math.}  pp. 1--47.
\bibitem{Oberman_08} A. Oberman (2008). Computing the convex envelope using a nonlinear partial differential equation, {\em Mathematical Models and Methods in Applied Sciences} (M3AS), Vol. 18. No 5 759--780.
\bibitem{Carlier_Galichon_11} G. Carlier and A. Galichon (2012). Exponential convergence for a convexifying equation, {\em ESAIM COCV}, 18, no. 3, 611--620.
\bibitem{KS} B. Kawohl and C. Schwab (1998). Convergent finite elements for a class of nonconvex variationnal problems, {\em IMA Journal of Numerical Analysis}  {\bf 18}:133-149.
\end{thebibliography}

\newpage

\newpage

\newpage

\begin{table}[htbp]
\begin{center}
\begin{tabular}{|c|c|c|c|}
\hline
\includegraphics{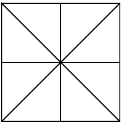} & \includegraphics{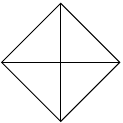} & \includegraphics{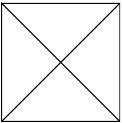} & \includegraphics{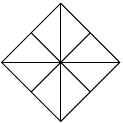}\\
\hline
$h^2\left(\begin{array}{cc} u_{xx} & 2u_{xy}\\ 2u_{xy} & u_{yy}\end{array}\right)$ &
$h^2\left(\begin{array}{cc} u_{xx} & 0\\ 0 & u_{yy}\end{array}\right)$ &
$4h^2/3\left(\begin{array}{cc} u_{xx} & u_{xy}\\ u_{xy} & u_{yy}\end{array}\right)$ &
$2h^2/3\left(\begin{array}{cc} u_{xx} & u_{xy}\\ u_{xy} & u_{yy}\end{array}\right)$ \\
\hline
\end{tabular}
\caption{First order of the Hessian for various meshes ($P_1-P_1$).}
\label{table1}
\end{center}
\end{table}


\end{document}